\documentclass[11pt]{article}

\def\mathbb{\Bbb}
\usepackage{amsfonts,latexsym,amstext}
\usepackage{amsmath}
\usepackage{amssymb}

\newtheorem{theorem}{Theorem}[section]
\newtheorem{lemma}[theorem]{Lemma}
\newtheorem{proposition}[theorem]{Proposition}

\newtheorem{hypothesis}[theorem]{Hypothesis}

\newtheorem{remark}[theorem]{Remark}

\makeatletter
\@addtoreset{equation}{section}

\makeatother

\def\qed{{\hfill\hbox{\enspace${ \square}$}} \smallskip}
\def\sqr#1#2{{\vcenter{\vbox{\hrule height .#2pt \hbox{\vrule
 width .#2pt height#1pt \kern#1pt \vrule
width .#2pt} \hrule height .#2pt}}}}
\def\square{\mathchoice\sqr54\sqr54\sqr{4.1}3\sqr{3.5}3}

\def\ds{\begin{displaystyle}}
\def\eds{\end{displaystyle}}
\def\dis{\displaystyle }
\def\<{\langle }
\def\>{\rangle }

\def\R{\mathbb R}

\def\E{\mathbb E}
\def\P{\mathbb P}

\def\calb{{\cal B}}

\def\cald{{\cal D}}

\def\calf{{\cal F}}

\def\calh{{\cal H}}

\def\calo{{\cal O}}
\def\calp{{\cal P}}

\def\call{{\cal L}}

\author{Marco Fuhrman
\\
Dipartimento di Matematica,
Politecnico di Milano\\
piazza Leonardo da Vinci 32, 20133 Milano, Italy\\
e-mail: marco.fuhrman@polimi.it
\\ \\
Ying Hu\\
IRMAR, Universit\'e Rennes 1\\ Campus de Beaulieu, 35042 Rennes
Cedex, France\\ e-mail: ying.hu@univ-rennes1.fr
\\ \\
Gianmario Tessitore
\\
Dipartimento di Matematica e Applicazioni, Universit\`a di Milano-Bicocca\\
Via Cozzi 53, 20125 Milano, Italy\\
e-mail: gianmario.tessitore@unimib.it }

\title{Stochastic maximum principle for optimal control of SPDEs}
\date{}
\author{Marco Fuhrman
\\
Dipartimento di Matematica,
Politecnico di Milano\\
via Bonardi 9, 20133 Milano, Italy\\
e-mail: marco.fuhrman@polimi.it
\\ \\
Ying Hu\\
IRMAR, Universit\'e Rennes 1\\ Campus de Beaulieu, 35042 Rennes
Cedex, France\\ e-mail: ying.hu@univ-rennes1.fr
\\ \\
Gianmario Tessitore
\\
Dipartimento di Matematica e Applicazioni, Universit\`a di Milano-Bicocca\\
Via Cozzi 53, 20125 Milano, Italy\\
e-mail: gianmario.tessitore@unimib.it }

\begin{document}

\maketitle

\begin{abstract}
We prove a version of the maximum principle, in the sense of
Pontryagin, for the optimal control of a stochastic partial differential equation
driven by a finite dimensional Wiener process.
The equation is formulated in a semi-abstract form that allows
direct applications to a large class of controlled stochastic parabolic equations.
We allow for a diffusion coefficient dependent on the control parameter,
and the space of control actions is general, so that in particular we need to introduce
two adjoint processes. The second
adjoint process takes values in a suitable space of operators on $L^4$.
\end{abstract}

\section{Introduction}
The  problem of finding sufficient conditions for optimality
for a stochastic optimal control problem with infinite dimensional state equation,
along the lines of the Pontryagin maximum principle,
was already addressed  in the early 80's in the pioneering paper  \cite{Be}.

Despite the fact that the finite dimensional analogue of the problem has been solved, in complete generality, more than 20 years ago
(see the well known paper by S. Peng
\cite{Pe}) the infinite dimensional case still  has  important open issues both on the side of the generality
of the abstract model and on the side of its
 applicability to systems modeled by stochastic partial differential equations (SPDEs).

 In particular, whereas the  Pontryagin maximum principle for infinite dimensional
 stochastic control problems is a well known result
 as far as the control domain is convex (or the diffusion does not depend on the control), see \cite{Be, HuPe2}, for the general case
(that is when the control domain need not be convex and the diffusion coefficient can contain a control variable)
existing results  are limited to abstract evolution equations under assumptions that are not satisfied
by the  large majority of concrete SPDEs.

The technical obstruction is related to the fact that (as it was pointed out in \cite{Pe}) if the control domain is not convex
the optimal control has to be perturbed by the so called ``spike variation".
Then if the control enters the diffusion,  the irregularity in time of the Brownian trajectories imposes to take into account
a second variation process. Thus   the stochastic maximum principle has to involve an adjoint process for the  second variation.
In the finite dimensional case such a process can be characterized as the solution of a matrix valued  backward stochastic differential equation
(BSDE) while in the infinite dimensional case the process naturally lives in a non-Hilbertian space of operators and its characterization
is much more difficult.
 Moreover  the applicability of the abstract results to concrete controlled SPDEs is another delicate step due to the specific
 difficulties that they involve such as the lack of regularity of
 Nemytskii-type coefficients in $L^p$ spaces.

 The present results (that were anticipated in the
 6th International Symposium on BSDEs and Applications  - Los Angeles, 2011 -
 and published in a short version in \cite{FuHuTe})
 are, as far as we know, the only ones that can cover, for instance, a controlled stochastic heat equation (with finite dimensional noise) such as:
\begin{equation}\label{state-equation}
    \left\{
\begin{array}{lll}
dX_t(x)&=&\dis A X_t(x)\,dt +b(x,X_t(x),u_t)\,dt+
\sum_{j=1}^m\sigma_j(x,X_t(x),u_t)\,d\beta^j_t,\qquad  t\in [0,T], x\in\calo,
\\
X_0(x)&=&x_0(x),
\end{array}
\right.
\end{equation}
with $A=\Delta_x$ with appropriate boundary conditions, and a cost functional as follows:
$$
J(u)=\E \int_0^T\int_\calo l(x,X_t(x),u_t)\,dx\,dt
+\E \int_\calo h(x,X_T(x))\,dx,
$$
 $\calo  \subset \R^n$ being a bounded open set with regular boundary.

 We stress the fact that in this paper the state equation is formulated, as above, only in a semi-abstract way  in order, on one side, to cope with all the difficulties
 carried by the concrete non-linearities and on the other to take advantage of the regularizing properties of the leading elliptic operator.

  Concerning other results on the infinite dimensional stochastic Pontryagin maximum principle, as we already mentioned in \cite{Be}   and \cite{HuPe2}
 the case of diffusion independent on the control is treated (with the difference that in \cite{HuPe2}
 a complete characterization of
 the adjoint process to the first variation as the unique mild solution to a suitable BSDE  is achieved).
 Then in \cite{Zh} the case of linear state equation and cost functional is addressed. In this case as
 well, the second variation process is not needed.

  The pioneering paper \cite{TaLi}
   is the first one in which the general case is addressed
  with, in addition, a general class of noises
 possibly with jumps. The adjoint process of the second variation  $(P_t)_{t\in [0,T]}$ is characterized as the solution of a
 BSDE in the (Hilbertian) space of Hilbert Schmidt  operators.
This forces  to assume a very strong regularity on the abstract
state equation and control functional that prevents application of the general results to SPDEs.
Recently in \cite{LuZh} $P_t$ was characterized as ``transposition solution" of a backward stochastic evolution equation
in $\mathcal{L}(L^2(\calo))$. Coefficients are required to be twice Fr\'echet-differentiable as operators in $L^2(\calo)$.
Finally even more recently in a couple of preprints
\cite{DuMe} \cite{DuMe2} the process $P_t$ is characterized in a similar way as it is in \cite{FuHuTe} and here. Roughly speaking it is characterized
as a
suitable stochastic bilinear form (see relation (\ref{defpt})).
As it is the case in \cite{LuZh}, in  \cite{DuMe} and \cite{DuMe2} as well the regularity assumptions on the coefficients are too restrictive
to apply directly  the general results to controlled SPDEs. On the other side in \cite{DuMe2} an unbounded diffusion term is included in the model
that can not be covered by the present results.
 Finally other variants of the problem have been studied. For instance in \cite{Gu}
a maximum principle for a SPDE
with  noise and control on the boundary but control independent diffusion is addressed,
see also \cite{OkSuZh} for a case with delay.

The paper is structured as follows. In Section 2 we fix notations and standing assumptions. In Section 3 we state the main result.
In Section 4 we recall the spike variation technique and introduce the first variation process, the corresponding adjoint process and the second variation process
together with crucial  estimates on them.  We stress here the fact that the structure of the second variation process forces to develop a theory in the
 $L^p$ spaces   for the state equation and its perturbations through spike variation.
In section 5 we complete the proof of the stochastic maximum principle. This is achieved by characterizing the adjoint of the second variation as a progressive process $(P_t)_{t\in [0,T]}$ with values in
 the space of linear
bounded operators $L^4\to (L^4)^*=L^{4/3}$. Namely $P_t$ is defined through the stochastic bilinear form
$$
    \< P_t f,g\>= \E^{\calf_t}\int_t^T\<\bar H_s Y_s^{t,f}, Y_s^{t,g}\>\, ds
+
 \E^{\calf_t}\< \bar
   h  Y_T^{t,f},    Y_T^{t,g}  \>, \qquad \P-a.s.
$$
where $( Y_s^{t,f})_{s\in [t,T]}$ is the mild solution of a suitable infinite dimensional forward stochastic equation
(see equation (\ref{auxiliaryy})).
The study of the regularity of  process $(P_t)_{t\in [0,T]}$ is one of the main technical issues of this paper (together with the $L^p$
estimates of the first and second variations)
and exploits the specific properties of the semigroup generated by the elliptic differential operator.
Finally in the Appendix we report some results on stochastic integration in $L^p$ spaces. For the reader's convenience we give complete and direct proofs of some results (including a version of the 
It\^{o}  inequality, see (\ref{isometryiito})). Such results are particular cases of the ones obtained in the framework of stochastic calculus in UMD Banach spaces, see \cite{vNVW}.

\section{Notations and preliminaries}

We begin by formulating an abstract form of the controlled PDE.

Let $(D,\cald,m)$ be
a measure space with finite measure
(in the applications
$D$ is an open subset of $\R^N$ and $m$ is the Lebesgue measure).
We will consider the usual real spaces
$L^p(D,\cald,m)$, $p\in [1,\infty)$, which are shortly denoted by
$L^p$ and endowed with the usual norm
$\|\cdot\|_p$.

Let $(W^1_t,\ldots,W^d_t)_{t\ge0}$ be a standard, $d$-dimensional
Wiener process defined in some complete probability space $(\Omega,\calf,\P)$.
We denote by $(\calf_t)_{t\ge0}$ the corresponding natural filtration,
augmented in the usual way, and we denote by $\calp$ the progressive $\sigma$-algebra
on $\Omega\times [0,\infty)$ (or on a finite interval $[0,T]$, by abuse of notation).
We will assume that there exist  regular conditional probabilities
$\P(\cdot|\calf_t)$
given any $\calf_t$: this holds for instance if the Wiener process is canonically
realized on the space of $\R^d$-valued continuous functions.

As the space of control actions we take  a separable metric space  $U$,
endowed with its Borel $\sigma$-algebra $\calb(U)$.
In general, we denote $\calb(\Lambda)$ the Borel $\sigma$-algebra
of any topological space $\Lambda$.
We fix a finite time horizon $T>0$
and by a control process we mean
any  progressive process $(u_t)_{t\in [0,T]}$ with values in $U$.

We consider the following controlled stochastic equation:
\begin{equation}\label{pdecontr}
    \left\{\begin{array}{l}
    dX_t(x)=AX_t(x)\,dt + b(t,x,X_t(x),u_t)\,dt
    + \sum_{j=1}^d \sigma_j(t,x,X_t(x),u_t)\,dW^j_t,
    \\
    X_0(x)=x_0(x)
\end{array}\right.
\end{equation}
and the cost functional
\begin{equation}\label{cost}
J(u)=\E \int_0^T\int_D
l(t,x,X_t(x),u_t)\,m(dx)\,dt+\E\,
\int_D h(x,X_T(x))\,m(dx).
\end{equation}
A control process $u$ is called optimal if it minimizes the cost over all
control processes. Denoting  by $X$  the corresponding trajectory
we say that $(u,X)$ is an optimal pair.

\begin{hypothesis}\label{standing}

\begin{enumerate}

\item The operator
$A$ is the infinitesimal generator of a strongly continuous semigroup $(e^{tA})_{t\ge 0}$
of linear bounded operators on $L^2$.
We assume that there exist constants $\bar p>8$ and $M>0$
such that for $p\in [2,\bar p]$ we have
$e^{tA}(L^p)\subset L^p$ and $\|e^{tA}f\|_p\le M\|f\|_p$
for every ${t\in[ 0,T]}$ and $f\in L^p$.

\item  For $\phi=b$ or $\phi=l$ or $\phi=\sigma_j$, $j=1,\ldots,d$, the functions
$$
\phi(\omega,t,x,r,u): \Omega\times [0,T]\times D\times \R \times U\to \R,
\quad
h(\omega,x,r): \Omega\times D\times \R \to \R,
$$
are assumed to be measurable with respect to $\calp\otimes \cald\otimes
\calb(\R) \otimes \calb(U)$ and $\calb(\R)$ (respectively,
$\calf_T\otimes \cald\otimes
\calb(\R) $ and $\calb(\R)$).
\item
For every $(\omega,t,x,u)$, the functions
$r\mapsto \phi(\omega,t,x,r,u)$ and
$r\mapsto h(\omega,x,r)$
are continuous and have first and second
derivatives, denoted $\phi'$ and $\phi''$ (respectively, $h'$ and $h''$), which are also continuous
functions of $r$.
We also assume that
$$
(|\phi'|+|\phi''|+|h'|+|h''|)
(\omega,t,x,r,u)
\le K,
$$
$$
(|\phi|+|h|)
(\omega,t,x,r,u)
\le K (|r|+|\bar\psi(x)|),
$$
for some constant $K$, some $\bar\psi\in L^{\bar p}$ and for all $(\omega,t,x,r,u)$.

\item $x_0\in L^{\bar p}$.

\end{enumerate}
\end{hypothesis}

From now on we adopt the convention of summation over repeated indices, so that
we will drop the symbol $\sum_{j=1}^d$ in (\ref{pdecontr}).

Under the stated assumptions, for every control process $u$ there exists
a unique solution of the state equation
(\ref{pdecontr}) in the so-called mild sense, i.e. an adapted
process $(X_t)_{t\in [0,T]}$ with values in $L^2$, with continuous
trajectories, satisfying $\P$-a.s.
\begin{equation}\label{pdecontrmild}
    X_t=e^{tA}x_0  + \int_0^t e^{(t-s)A}b(s,\cdot,X_s(\cdot),u_s)\,ds
    + \int_0^t e^{(t-s)A} \sigma_j(s,\cdot ,X_s(\cdot),u_s)\,dW^j_s,
    \qquad t\in [0,T].
\end{equation}
Here and below, equalities like \eqref{pdecontrmild} are understood
to hold $m$-a.e., and uniqueness is understood up to modification
of $L^2$-valued random processes.

We need to prove the following higher summability
property of trajectories.
We introduce the notation
\begin{equation}\label{normpprocess}
|||X|||_p=
    \sup_{t\in[0,T]}(\E\|X_t\|_p^p)^{1/p}.
\end{equation}

\begin{proposition}\label{lptrajectory}
For every $p\in [2,\bar p]$, $(X_t)_{t\in [0,T]}$ is
a progressive
process with values in $L^p$, satisfying
$
\sup_{t\in[0,T]} \E\|X_t\|_p^p<\infty.
$
\end{proposition}

\noindent{\bf Proof.}
 We consider the Banach space
  of progressive $L^p$-valued
processes $(X_t)_{t\in [0,T]}$ such that the norm
$|||X|||_p$ is finite.
    For such a process $X$ we define
    $$
\Gamma(X)_t= e^{tA}x_0+
  \int_0^t e^{(t-s)A}
  b(s,\cdot ,X_s(\cdot),u_s)\,ds
    + \int_0^t e^{(t-s)A}
       \sigma_j(s,\cdot ,X_s(\cdot),u_s)
\,dW^j_s.
$$
It can be proved that the map $\Gamma$ is a contraction with
respect to the  norm $|||\cdot|||_p$, provided
$T$ is sufficiently small. Its unique fixed point is then the required solution.
    The restriction on $T$ is then removed in a standard way by subdividing $[0,T]$
    into appropriate subintervals.

The fact that $\Gamma$ is a well defined  contraction
follows from moment estimates of the
stochastic integrals in $L^p$.  We limit ourselves to showing the contraction property
assuming for simplicity $b=0$. In this case, if $|||X|||_p+|||Y|||_p<\infty$, we have
by \eqref{isometryiitoter}
$$
\E\|
\Gamma(X)_t-\Gamma(Y)_t\|_p^p\le
 c_p\int_0^t \E \|
 e^{(t-s)A}[\sigma(s,\cdot ,X_s(\cdot),u_s)
 -\sigma(s,\cdot ,Y_s(\cdot),u_s)
 ]
 \|^p_{L^p(D;\R^d)}ds\; t^{(p-2)/2}.
 $$
Using the $L^p$-boundedness
of $e^{tA}$ and the Lipschitz character of $\sigma$ which follows from Hypothesis
\ref{standing}-$3$ we obtain
$$
\E\|
\Gamma(X)_t-\Gamma(Y)_t\|_p^p\le C
  \int_0^t \E \|
X_s-Y_s \|^pds\; t^{(p-2)/2} \le C
|||X-Y|||_p^p T^{p/2}
 $$
for some constant $C$ independent of $T$.
The contraction property
follows immediately for $T$ sufficiently small.
\qed

\section{Statement of the main result}

\subsection{Statement of the stochastic maximum principle}

For our main result  we also need the following assumptions.

\begin{hypothesis}\label{basisinl4}
There exists a complete orthonormal basis $(e_i)_{i\ge 1}$ in
$L^2$ which is also a Schauder basis of $L^4$.
\end{hypothesis}

\begin{hypothesis}\label{regularityinl4}
The restriction of
$(e^{tA})_{t\ge 0}$ to the space $L^4$ is a strongly continuous
analytic semigroup and the domain of its
 infinitesimal generator is compactly embedded in $L^4$.
\end{hypothesis}

We note that Hypothesis \ref{basisinl4} is satisfied
for a large class of
measure spaces $(D,\cald,m)$, typically with a basis of Haar type.

In the following a basic role will be played by the space
of linear
bounded operators $L^4\to (L^4)^*=L^{4/3}$ endowed with the usual operator norm,
that we simply denote by $\call$. Clearly, $\call$ may be identified with the
space of bounded bilinear forms on $L^4$.
The duality between $g\in L^4$ and $h\in L^{4/3}$
will be denoted $\<h,g\>$. As it is customary when dealing with
spaces of operators endowed with the operator norm,
when considering random variables or processes with values
in $\call$, the latter will be endowed with the Borel $\sigma$-algebra of the weak topology
(the weakest topology making all the functions $T\mapsto \<Tf,g\>$ continuous, $f,g\in L^4$);
note that this  is in general different from the Borel $\sigma$-algebra of the
topology corresponding to the operator norm $\|T\|_\call$.

For $u\in U$ and $X,p,q^1,\ldots,q^d\in L^2$  denote
$$
\calh(t,u,X,p,q^1,\ldots,q^d)=
\int_D [ l(t,x,X(x),u)+
 b(t,x,X(x),u)
p(x)
+
\sigma_{j}(t,x,X(x),u)
q^{j}(x)]
\,m(dx)
$$

\begin{theorem}\label{mainresult}
Let
$({X},{u})$ be an optimal pair. Then there exist
 progressive
processes
$(P_t)_{t\in [0,T]}$ and $(p_t,q_t^1,\ldots, q_t^d )_{t\in [0,T]}$,
with values in $\call$ and $(L^2)^{d+1}$  respectively,
for which the following inequality holds, $\P$-a.s. for a.e. $t\in
[0,T]$:
 for every $v\in U$,
$$
\begin{array}{l}
\calh(t,v,{X}_t,p_t, q_{t}^1,\ldots,q_{ t}^d)-
 \calh(t, {u}_t, {X}_t,p_t, q_{t}^1,\ldots,q_{t}^d)
 \\\dis
 \qquad
+\frac{1}{2}\,\<P_t[\sigma_j(t,\cdot, {X}_t(\cdot),v)-
\sigma_j(t,\cdot, {X}_t(\cdot), {u}_t)],
\sigma_j(t,\cdot, {X}_t(\cdot),v)-
\sigma_j(t,\cdot, {X}_t(\cdot), {u}_t)\>
\ge 0.
\end{array}
$$

The process $(p,q^1,\ldots, q^d )$
satisfies
$
\sup_{t\in[0,T]}\E\|p_t\|_2^2+
\E\int_0^T\sum_{j=1}^d\|q^j_t\|_2^2\,dt
<\infty,
$
and
it is the unique solution to equation
\eqref{firstadjointmild} below.

The  process $P$ satisfies
$    \sup_{t\in [0,T]}\E\|P_t\|_\call^2<\infty$ and
it is defined in
Proposition \ref{processp} below (formula \eqref{defpt}).

 \end{theorem}

$(p,q^1,\ldots, q^d )$ and $P$ will be called
the first and second adjoint process, respectively.

\subsection{Application to stochastic PDEs of parabolic type}

The purpose of this short subsection is to show
that  the main result can be immediately applied
to concrete cases of controlled stochastic
PDE of parabolic type on domains of Euclidean space.

Let $D$ be a bounded open subset of $\R^n$ with smooth
boundary $\partial D$, and let $m$ be the Lebesgue
measure. Consider the following
PDE of reaction-diffusion
type
\begin{equation}\label{pdecontrconcrete}
    \left\{\begin{array}{l}
    dX_t(x)=\Delta X_t(x)\,dt + b(t,x,X_t(x),u_t)\,dt
    + \sum_{j=1}^d \sigma_j(t,x,X_t(x),u_t)\,dW^j_t, \quad t\in[0,T], x\in D,
    \\
    X_t(x)=0,\quad t\in[0,T], x\in \partial D,
    \\
    X_0(x)=x_0(x),\quad x\in D,
\end{array}\right.
\end{equation}
and the cost functional
\eqref{cost}.
In this example
 the Wiener process, the space of control actions
and the space of  control processes are as  before; on the coefficients
$b,\sigma_j,l,h$ we make the assumptions of
Hypothesis \ref{standing}, points 2 and 3; finally we suppose $x_0\in L^{\bar p}$
for some $\bar p>8$.

We claim that all the conclusions of Theorem \ref{mainresult} hold true.

Indeed, we can define the operator $A=\Delta$ as an unbounded operator in $L^2$
with domain $H^2(D)\cap H^1_0(D)$ (the standard Sobolev spaces).
Then $A$ generates a strongly continuous, analytic contraction semigroup in
all the spaces $L^p$, $1<p<\infty$, and the domain of $A$ is compactly
embedded: see e.g. \cite{lu} or \cite{pa}. Therefore  Hypothesis \ref{standing}, point 1,
and Hypothesis \ref{regularityinl4} hold true. Finally, Hypothesis \ref{basisinl4}
is verified by a Haar basis.

By similar arguments
instead of $\Delta$ one can consider more general  operators of elliptic type
with appropriate boundary conditions.

\section{The spike  variation method and the first adjoint process}

Throughout this section we assume that Hypothesis \ref{standing} holds,
whereas Hypotheses \ref{basisinl4} and \ref{regularityinl4} will be needed only starting
from the next section.

\subsection{Spike  variation method and expansion of the state and the cost}

Suppose that $u$ is an optimal control and $X$ the corresponding optimal
trajectory.  We fix  $t_0\in (0,T)$ and $\epsilon>0$ such that
$[t_0,t_0+\epsilon]\subset (0,T)$, we fix a control process $v$ and
we introduce in the usual way the spike variation process
$$
u^\epsilon_t=\left\{\begin{array}{ll}
v_t, & {\rm if} \; t\in [t_0,t_0+\epsilon],
\\
u_t, & {\rm if} \; t\notin [t_0,t_0+\epsilon].
\end{array}\right.
$$
We denote by $X^\epsilon$ the trajectory corresponding to $u^\epsilon$.
We are going to construct two $L^2$-valued stochastic processes,
denoted $Y^\epsilon$ and $Z^\epsilon$, in such a way that the difference
$X-X^\epsilon-Y^\epsilon-Z^\epsilon$ is small (in the sense of
Proposition \ref{variations} below) and the difference of the cost
functional $J(u^\epsilon)-J(u)$ can be expressed in an appropriate
form involving $Y^\epsilon$ and $Z^\epsilon$
up to a small remainder: see
Proposition \ref{costvariationprelim}.

Define
$$
\begin{array}{lll}
\delta^\epsilon \sigma_j (t,x)&=&
\sigma_j(t,x,X_t(x),u^\epsilon_t)- \sigma_j(t,x,X_t(x),u_t)
\end{array}
$$
and consider the stochastic PDE
\begin{equation}\label{firstvariation}
    \left\{\begin{array}{lll}
    dY^\epsilon_t(x)&=&\dis
    AY^\epsilon_t(x)\,dt + b'(t,x,X_t(x),u_t)Y^\epsilon_t(x)\,dt
    \\
 &&  \dis +  \sigma_j'(t,x,X_t(x),u_t)Y^\epsilon_t(x)\,dW^j_t
    + \delta^\epsilon \sigma_j (t,x)\,dW^j_t,
    \\
    Y^\epsilon_0(x)&=&0.
\end{array}\right.
\end{equation}
By the standard theory of stochastic evolution equations in Hilbert
spaces, see e.g. \cite{DaZa},
there exists
a unique solution to
(\ref{firstvariation}) in the  mild sense, i.e. a progressive
process $(Y^\epsilon_t)_{t\in [0,T]}$ with values in $L^2$, satisfying
$\sup_{t\in[0,T]}\E\|Y^\epsilon_t\|_2^2<\infty$
and, for every $t\in[0,T]$,
\begin{equation}\label{firstvariationmild}
\begin{array}{lll}
    Y^\epsilon_t&=&\dis
  \int_0^t e^{(t-s)A}[
    b'(s,\cdot,X_s(\cdot),u_s)Y^\epsilon_s(\cdot)
    + \delta^\epsilon b(s,\cdot)]\,ds
    \\
&&\dis
    + \int_0^t e^{(t-s)A}
    [ \sigma_j'(s,\cdot,X_s(\cdot),u_s)Y^\epsilon_s(\cdot)
    + \delta^\epsilon \sigma_j (s,\cdot)]
\,dW^j_s,
    \qquad \P-a.s.
    \end{array}
\end{equation}
For the sequel we need the following more precise result.

\begin{proposition}\label{lpfirstvariation}
For every $p\in [2,\bar p]$, $(Y^\epsilon_t)_{t\in [0,T]}$ is
a progressive
process with values in $L^p$, satisfying
$$
|||Y^\epsilon|||_p=
\sup_{t\in[0,T]}(\E\|Y^\epsilon_t\|_p^p)^{1/p}\le C\epsilon^{1/2}.
$$
\end{proposition}

To prove this result we need the following lemma, that will be used
several times.
\begin{lemma}\label{lineqereqestimates}
Given $\calp\otimes \cald$-measurable processes
$\bar  a, \bar \alpha, \bar b^j , \bar \beta^j $ consider the linear equation:
\begin{equation}\label{lineqereq}
    \left\{\begin{array}{lll}
    dV _t(x)&=&\dis
    AV _t(x)\,dt + \bar a  (t,x)V _t(x)\,dt +\bar \alpha  (t,x)\,dt
      +  \bar b^j (t,x)V _t(x)\,dW^j_t
    + \bar \beta^j (t,x)\,dW^j_t,
    \\
V_0(x)&=&0.
\end{array}\right.
\end{equation}
Suppose $\bar  a ,\bar b^j$ bounded and
 $p\in[2,\bar p]$. Then the following holds.
 \begin{enumerate}
   \item There exists
a unique solution to \eqref{lineqereq}
in the  mild sense, i.e. a progressive
process $(V_t)_{t\in [0,T]}$ with values in $L^p$, satisfying
\begin{equation}\label{estimatelineareq}
|||V|||_p=
    \sup_{t\in[0,T]}(\E\|V_t\|_p^p)^{1/p}\le
    C\,
     \int_0^T (\E\|\bar \alpha _t\|_p^p)^{1/p}dt
    +
    C\,
    \left(\int_0^T (\E\|\bar \beta_t\|_p^p)^{2/p}\right)^{1/2}
\end{equation}
and, for every $t\in[0,T]$,
$$
    V_t=
  \int_0^t e^{(t-s)A}[
    \bar a  (s,\cdot )V_s(\cdot)
    + \bar \alpha  (s,\cdot)]\,ds
    + \int_0^t e^{(t-s)A}
[      \bar b _j (s,\cdot)V_s(\cdot)
    + \bar \beta^j (s,\cdot)]
\,dW^j_s,
    \qquad \P-a.s.
    $$
provided the right-hand side of \eqref{estimatelineareq} is finite.

   \item If in addition     $\sup_{t\in[0,T]}\E(\|\bar \alpha _t\|_p^p+
\|\bar \beta_t\|_p^p) <\infty$ and
$\bar  a ,\bar b^j$ are supported in a time interval of length $\epsilon$, then
$$
    \sup_{t\in[0,T]}(\E\|V_t\|_p^p)^{1/p}\le
    C\, \epsilon\,
    \sup_{t\in[0,T]}(\E\|\bar \alpha_t\|_p^p)^{1/p}
    +
    C\,\sqrt{\epsilon}\,
 \sup_{t\in[0,T]}(\E\|\bar \beta_t\|_p^p)^{1/p}.
$$
or equivalently
\begin{equation}\label{estimatelineareqepsilon}
|||V|||_p\le
    C\, \epsilon\,
    |||\bar \alpha|||_p
    +
    C\,\sqrt{\epsilon}\,
|||\bar \beta|||_p.
\end{equation}

   \item In the case $p=2$ we have
   \begin{equation}\label{estimatelineareq2}
    \sup_{t\in[0,T]}\E\|V_t\|_2^2 \le
    C\,
     \int_0^T  \E\|\bar \alpha _t\|_2^2 \,dt
    +
    C\,
 \int_0^T  \E\|\bar \beta_t\|_2^2 dt
 =C\, (\|\bar \alpha \|_{L^2(\Omega\times D\times [0,T])}^2
    +\|\bar \beta\|_{L^2(\Omega\times D\times [0,T])}^2).
\end{equation}

 \end{enumerate}

In \eqref{estimatelineareq}, \eqref{estimatelineareqepsilon}, \eqref{estimatelineareq2}
we set  $\|\bar \beta_t\|_p:=\|\bar \beta_t\|_{L^p(D;\R^d)}$, and
 the constant $C$ depends on the bounds
on $\bar  a ,\bar b^j$, on the semigroup $(e^{tA})$, and on $p$ and $T$.
\end{lemma}

 \noindent {\bf Proof.} We consider again the Banach space
  of progressive $L^p$-valued
processes $(V_t)_{t\in [0,T]}$ endowed with the norm
$|||V|||_p=
    \sup_{t\in[0,T]}(\E\|V_t\|_p^p)^{1/p}$.
By the same arguments as in the proof of
Proposition \ref{lptrajectory}  we can prove that
the map $\Gamma$ defined as
    $$
\Gamma(V)_t=
  \int_0^t e^{(t-s)A}
    \bar a  (s,\cdot )V_s(\cdot)\,ds
    + \int_0^t e^{(t-s)A}
       \bar b _j (s,\cdot)V_s(\cdot)
\,dW^j_s
$$
is a (linear) contraction with respect to the norm $|||\cdot|||_p$, provided
$T$ is sufficiently small.
    Therefore there exists a unique solution $V$ and it satisfies the inequality
    $$
    |||V|||_p\le C\, |||
  \int_0^\cdot e^{(\cdot-s)A}
    \bar \alpha  (s)\,ds|||_p+ |||
     \int_0^\cdot e^{(\cdot-s)A}
     \bar \beta^j (s)
\,dW^j_s|||_p.
    $$
The inequality \eqref{estimatelineareq} follows from an estimate of those
stochastic integrals, using \eqref{isometryiito} and the $L^p$-boundedness
of $e^{tA}$.
    The restriction on $T$ is then removed by subdividing $[0,T]$
    into appropriate subintervals.
 Finally,
 \eqref{estimatelineareqepsilon} and  \eqref{estimatelineareq2} follow from
 \eqref{estimatelineareq} and the H\"older inequality.
 \qed

 \noindent {\bf Proof of Proposition  \ref{lpfirstvariation}.}
 This is an immediate corollary of the previous lemma, noting that
 $|||\delta^\epsilon \sigma_j |||_p\le C$ as a consequence of the linear growth
 condition on $\sigma_j $
 (Hypothesis \ref{standing}-$4$) and the fact that $|||X|||_p<\infty$ by Proposition
\ref{lptrajectory}.
 \qed

Define
$$
\begin{array}{lll}
\delta^\epsilon b(t,x)&=&
b(t,x,X_t(x),u^\epsilon_t)- b(t,x,X_t(x),u_t),
\\
\delta^\epsilon b'(t,x)&=&
b'(t,x,X_t(x),u^\epsilon_t)- b'(t,x,X_t(x),u_t),
\\
\delta^\epsilon \sigma_j' (t,x)&=&
\sigma_j'(t,x,X_t(x),u^\epsilon_t)- \sigma_j'(t,x,X_t(x),u_t),
\end{array}
$$
and consider the following stochastic PDE:
\begin{equation}\label{secondvariation}
    \left\{\begin{array}{lll}
    dZ^\epsilon_t(x)&=&\dis
    AZ^\epsilon_t(x)\,dt + b'(t,x,X_t(x),u_t)Z^\epsilon_t(x)\,dt
    \\
 &&  \dis
    +\frac12 b''(t,x,X_t(x),u_t)Y^\epsilon_t(x)^2\,dt
    +     \delta^\epsilon b(t,x)\,dt
    + \delta^\epsilon b'(t,x)Y^\epsilon_t(x)\,dt
    \\
 &&  \dis +  \sigma_j'(t,x,X_t(x),u_t)Z^\epsilon_t(x)\,dW^j_t
       +\frac12 \sigma_j''(t,x,X_t(x),u_t)Y^\epsilon_t(x)^2\,dW^j_t
       \\
 &&  \dis
    + \delta^\epsilon \sigma_j' (t,x)Y^\epsilon_t(x)\,dW^j_t,
    \\
    Z^\epsilon_0(x)&=&0
\end{array}\right.
\end{equation}
By the standard theory
there exists
a unique solution to
(\ref{secondvariation}) in the  mild sense, i.e. a progressive
process $(Z^\epsilon_t)_{t\in [0,T]}$ with values in $L^2$, satisfying
$\sup_{t\in[0,T]}\E\|Z^\epsilon_t\|_2^2<\infty$
and, for every $t\in[0,T]$,
\begin{equation}\label{secondvariationmild}
\begin{array}{lll}
    Z^\epsilon_t&=&\dis
  \int_0^t e^{(t-s)A}[
    b'(s,\cdot,X_s(\cdot),u_s)Z^\epsilon_s(\cdot)
    +\frac12 b''(s,\cdot,X_s(\cdot),u_s)Y^\epsilon_s(\cdot)^2
    + \delta^\epsilon b(s,\cdot)+ \delta^\epsilon b'(s,\cdot)Y^\epsilon_s(\cdot)]\,ds
    \\
&&\dis
    + \int_0^t e^{(t-s)A}
    [ \sigma_j'(s,\cdot,X_s(\cdot),u_s)Z^\epsilon_s(\cdot)
     +\frac12 \sigma_j''(s,\cdot,X_s(\cdot),u_s)Y^\epsilon_s(\cdot)^2
    + \delta^\epsilon \sigma_j' (s,\cdot)Y^\epsilon_s(\cdot)]
\,dW^j_s,
    \, \P\text{\it-a.s.}
    \end{array}
\end{equation}

For the sequel we need the following result.

\begin{proposition}\label{lpsecondvariation}
For every $p\in [2,\bar p /2]$, $(Z^\epsilon_t)_{t\in [0,T]}$ is
a progressive
process with values in $L^p$, satisfying
$$
|||Z^\epsilon|||_p=
\sup_{t\in[0,T]}(\E\|Z^\epsilon_t\|_p^p)^{1/p}\le C\epsilon.
$$
\end{proposition}

\noindent{\bf Proof.}
The result follows from
Lemma \ref{lineqereqestimates} applied to equation \eqref{secondvariation}.
In particular inequality
\eqref{estimatelineareq} shows that
$$
\begin{array}{lll}\dis
|||Z^\epsilon|||_p&\le&\dis
   C  \int_0^T (\E\|
     \frac12 b''( t,X_t ,u_t )(Y_t^\epsilon )^2
    + \delta^\epsilon b(t)+ \delta^\epsilon b'(t)Y_t^\epsilon
     \|_p^p)^{1/p}dt
    \\&&\dis
    +
    C
    \left(\int_0^T (\E\|
    \frac12 \sigma''(t, X_t ,u_t )(Y_t^\epsilon)^2
    + \delta^\epsilon \sigma'(t)  Y_t^\epsilon
    \|_p^p)^{2/p}\right)^{1/2}.
\end{array}
$$
The proof is now concluded  estimating the right-hand side of this inequality.

Since
$
\| \sigma''(t, X_t ,u_t )(Y_t^\epsilon)^2
    \|_p \le C \|
(Y_t^\epsilon)^2
    \|_p= C \|
Y_t^\epsilon
    \|_{2p}^2$
we have
$$
    \left(\int_0^T (\E\|
 \sigma''(t, X_t ,u_t )(Y_t^\epsilon)^2
    \|_p^p)^{2/p}\right)^{1/2}
     \le
    C
\left(\int_0^T (\E\|
 Y_t^\epsilon
    \|_{2p}^{2p})^{2/p}\right)^{1/2}
    \le C|||Y^\epsilon|||_{2p}^2\le C\epsilon,
$$
by Proposition \ref{lpfirstvariation}, since
  $2p\le \bar p$.

Next we note that
 $|||\delta^\epsilon b |||_p\le C$, as a consequence of the linear growth
 condition on $b$
 (Hypothesis \ref{standing}-$3$) and the fact that $|||X|||_p<\infty$ by Proposition
\ref{lptrajectory}. Since $\delta^\epsilon b$ is supported in $[t_0,t_0+\epsilon]$
it follows that
$ \int_0^T (\E\|
     \delta^\epsilon b(t)
     \|_p^p)^{1/p}dt\le C\epsilon $.

The other terms are treated in a similar way.
\qed

\begin{proposition}\label{variations} We have
$$
\sup_{t\in[0,T]}(\E\|X^\epsilon_t-X_t-Y^\epsilon_t-Z^\epsilon_t\|_2^2)^{1/2}=o(\epsilon).
$$
\end{proposition}

As usual, $o(\epsilon)$ denotes any function of $\epsilon$ such that
$o(\epsilon)/\epsilon \to 0$ as $\epsilon\to 0$.
During the proof we will use the Taylor formula in the following form:
for a twice continuously differentiable real function $g$ on $\R$, and for
$r,h\in\R$,
\begin{equation}\label{taylor}
    g(r+h)=g(r)+g'(r)h+\int_0^1\int_0^1 g''(r+\lambda\mu h)\,\mu d\lambda d\mu\, h^2.
\end{equation}
Since  $\int_0^1\int_0^1 \mu d\lambda d\mu=1/2$
this can also be written
\begin{equation}\label{taylordue}
    g(r+h)=g(r)+g'(r)h+\frac12 g''(r) h^2+
    \int_0^1\int_0^1 [g''(r+\lambda\mu h)-g''(r)]\,\mu d\lambda d\mu\, h^2.
\end{equation}

{\bf Proof.}
We set $R^\epsilon= Y^\epsilon+Z^\epsilon$. We first show that $X+R^\epsilon$
is a solution of  the following equation in $L^2$:
\begin{equation}\label{eqxplusr}
\begin{array}{l}\dis
 X_t+R^\epsilon_t=e^{tA}x_0
 + \int_0^t e^{(t-s)A}b(s,\cdot,X_s(\cdot)+R^\epsilon_s(\cdot),u_s)\,ds
 -  \int_0^t e^{(t-s)A}G^\epsilon(s,\cdot)\,ds
 \\\dis
    + \int_0^t e^{(t-s)A} \sigma_j(s,\cdot ,X_s(\cdot)+R^\epsilon_s(\cdot),u_s)\,dW^j_s
   -  \int_0^t e^{(t-s)A}\Lambda_j^\epsilon(s,\cdot)\,dW^j_s,
\end{array}
\end{equation}
where $G^\epsilon=G ^{\epsilon,1}+G ^{\epsilon,2}+G ^{\epsilon,3}$,
 $\Lambda_j^\epsilon=\Lambda_j^{\epsilon,1}+\Lambda_j^{\epsilon,2}
 +\Lambda_j^{\epsilon,3}$,
$$
\begin{array}{l}\dis
G^{\epsilon,1}(s,x)= \int_0^1\int_0^1 [
b''(s,x,X_s(x)+\lambda \mu R^\epsilon_s(x),u_s^\epsilon)- b''(s,x,X_s(x),u_s)
] \,\mu d\lambda d\mu\,R^\epsilon_s(x)^2,
\\\dis
G ^{\epsilon,2}(s,x)= \frac12 b''(s,x,X_s(x),u_s)
\,(Z^\epsilon_s(x)^2+2Y^\epsilon_s(x)Z^\epsilon_s(x)),
\qquad
G ^{\epsilon,3}(s,x)= \delta^\epsilon b'(s,x)Z^\epsilon_s(x),
\\\dis
\Lambda_j^{\epsilon,1}(s,x)= \int_0^1\int_0^1 [
\sigma_j''(s,x,X_s(x)+\lambda \mu R^\epsilon_s(x),u^\epsilon_s)
- \sigma_j''(s,x,X_s(x),u_s)
] \,\mu d\lambda d\mu\,R^\epsilon_s(x)^2,
\\\dis
\Lambda_j ^{\epsilon,2}(s,x)= \frac12 \sigma_j''(s,x,X_s(x),u_s)
\,(Z^\epsilon_s(x)^2+2Y^\epsilon_s(x)Z^\epsilon_s(x)),
\qquad
\Lambda_j ^{\epsilon,3}(s,x)= \delta^\epsilon \sigma_j'(s,x)Z^\epsilon_s(x).
\end{array}
$$
To verify \eqref{eqxplusr} we use the Taylor formula \eqref{taylor} and obtain
\begin{eqnarray}
  b(s,x,X_s(x)+R^\epsilon_s(x),u_s^\epsilon) &\!=\!&
   b(s,x,X_s(x),u_s^\epsilon)+
   b'(s,x,X_s(x),u_s^\epsilon)R^\epsilon_s(x)\nonumber
   \\
    &&+ \int_0^1\int_0^1
b''(s,x,X_s(x)+\lambda \mu R^\epsilon_s(x),u_s^\epsilon)
 \,\mu d\lambda d\mu\,R^\epsilon_s(x)^2,
 \label{variaauxuno}
\\
  \sigma_j(s,x,X_s(x)+R^\epsilon_s(x),u_s^\epsilon)
  &\!=\!& \sigma_j(s,x,X_s(x),u_s^\epsilon)+
   \sigma_j'(s,x,X_s(x),u_s^\epsilon)R^\epsilon_s(x)\nonumber
   \\
    &&+ \int_0^1\int_0^1
\sigma_j''(s,x,X_s(x)+\lambda \mu R^\epsilon_s(x),u_s^\epsilon)
 \,\mu d\lambda d\mu\,R^\epsilon_s(x)^2.
 \label{variaauxdue}
\end{eqnarray}
We apply  $e^{(t-s)A}$ to \eqref{variaauxuno} and integrate $\int_0^t \,ds$,
we apply  $e^{(t-s)A}$ to \eqref{variaauxdue} and integrate $\int_0^t \,dW^j_s$,
and we add the resulting equalities. Comparing with \eqref{pdecontrmild},
\eqref{firstvariationmild} and
\eqref{secondvariationmild} we obtain \eqref{eqxplusr}.

Since $X^\epsilon$ is the trajectory corresponding to $u^\epsilon$ we have
$$
 X_t^\epsilon=e^{tA}x_0  +
 \int_0^t e^{(t-s)A}b(s,\cdot,X_s^\epsilon(\cdot),u^\epsilon_s)\,ds
    + \int_0^t e^{(t-s)A} \sigma_j(s,\cdot ,X_s^\epsilon(\cdot),u^\epsilon_s)\,dW^j_s.
    $$
Comparing with \eqref{eqxplusr} we see that $\Delta^\epsilon:= X^\epsilon - X-
R^\epsilon$ solves
$$
\begin{array}{l}\dis
 \Delta^\epsilon_t=
  \int_0^t e^{(t-s)A}\bar b^\epsilon(s,\cdot)\Delta^\epsilon_s(\cdot)\,ds
 +  \int_0^t e^{(t-s)A}G^\epsilon(s,\cdot)\,ds
 \\\dis
 \qquad   + \int_0^t e^{(t-s)A} \bar \sigma^\epsilon_j(s,\cdot )\Delta^\epsilon_s(\cdot)\,dW^j_s
   +  \int_0^t e^{(t-s)A}\Lambda_j^\epsilon(s,\cdot)\,dW^j_s,
\end{array}
$$
where
$$
\begin{array}{l}\dis
\bar b^\epsilon(s,x)=
  \int_0^1
 b'(s,x,X_s(x)+R^\epsilon_s(x)+\lambda  \Delta^\epsilon_s(x),u_s^\epsilon)
  \,d\lambda,
 \\\dis
 \bar \sigma_j^\epsilon(s,x)=
  \int_0^1
 \sigma_j'(s,x,X_s(x)+R^\epsilon_s(x)+\lambda  \Delta^\epsilon_s(x),u_s^\epsilon)
  \,d\lambda,
 \end{array}
$$
are bounded coefficients, uniformly in $\epsilon$.
We can then apply Lemma \ref{lineqereqestimates} and
specifically inequality \eqref{estimatelineareq2} arriving at
$$
\sup_{t\in[0,T]}(\E\|X^\epsilon_t-X_t-Y^\epsilon_t-Z^\epsilon_t\|_2^2)^{1/2}=
    \sup_{t\in[0,T]}\E\|\Delta^\epsilon_t\|_2^2 \le
    C\, (\|G^\epsilon\|_{L^2(\Omega\times D\times [0,T])}^2
    +\|\Lambda^\epsilon\|_{L^2(\Omega\times D\times [0,T])}^2).
$$
To finish the proof
 it remains to verify that the $L^2(\Omega\times D\times [0,T])$-norm
 of each term $G ^{\epsilon,i}$,
 $\Lambda_j^{\epsilon,i}$ ($i=1,2,3$) is $o(\epsilon)$.

 Let us verify that $\|G^{\epsilon,1}\|_{L^2(\Omega\times D\times [0,T])}=o(\epsilon)$.
 Write $G^{\epsilon,1}(t,x)=Q^{\epsilon}_t(x)R^{\epsilon}_t(x)^2$ where
 $$
 Q^{\epsilon}_t(x)=\int_0^1\int_0^1 [
b''(t,x,X_t(x)+\lambda \mu R^\epsilon_t(x),u^\epsilon_t)- b''(t,x,X_t(x),u_t)
] \,\mu d\lambda d\mu.
 $$
 Next take $p\in (2, \bar p/4]$, which is possible because $\bar p>8$, and let $q>1$
 be such that $\frac12=\frac1p+\frac1q$. Then
 $$
 \|G^{\epsilon,1}\|_{L^2(\Omega\times D\times [0,T])}\le
 \|Q^{\epsilon}\|_{L^q(\Omega\times D\times [0,T])}
 \|(R^{\epsilon})^2\|_{L^p(\Omega\times D\times [0,T])}.
 $$
 Since
 $$
 \|(R^{\epsilon})^2\|_{L^p(\Omega\times D\times [0,T])}=
 \|R^{\epsilon}\|^2_{L^{2p}(\Omega\times D\times [0,T])}
 \le C |||R^{\epsilon}|||^2_{ 2p}
 \le C (|||Y^{\epsilon}|||^2_{ 2p}+|||Z^{\epsilon}|||^2_{ 2p})
 \le C (\epsilon+\epsilon^2)
 $$
 by Propositions \ref{lpfirstvariation} and \ref{lpsecondvariation},
 it remains to show that $\|Q^{\epsilon}\|_{L^q(\Omega\times D\times [0,T])}\to 0$
 as $\epsilon\to 0$. We argue by contradiction: assume that there exists
 $\delta >0$ and a sequence $\epsilon_n\to 0$ such that
 \begin{equation}\label{contrconvdom}
  \E\int_0^T\int_D\left|
 \int_0^1\int_0^1 [
b''(t,x,X_t(x)+\lambda \mu R^\epsilon_t(x),u^\epsilon_t)- b''(t,x,X_t(x),u_t)
] \,\mu d\lambda d\mu\right|^qm(dx)dt
 \ge \delta.
 \end{equation}
 Since
 $
 \|R^{\epsilon}\|_{L^2(\Omega\times D\times [0,T])}
 \le C |||R^{\epsilon}|||_{ 2}
 \le C (|||Y^{\epsilon}|||_{ 2}+|||Z^{\epsilon}|||_{ 2})\to 0,
 $
there exists a subsequence $\epsilon_{n_k}$ such that
 $R^{\epsilon_{n_k}}\to 0$ a.s. with respect to the product measure
 $\P(d\omega)m(dx)dt$. Since $r\mapsto b''(t,x,r,u)$ is continuous,
 and due to the special definition of $u^\epsilon$, it follows that
 $b''(t,x,X_t(x)+\lambda \mu R^{\epsilon_{n_k}}_t(x),
 u^{\epsilon_{n_k}}_t)\to b''(t,x,X_t(x),u_t)$ a.s. with respect to

 $\P(d\omega)m(dx)dtd\lambda d\mu$. By dominated convergence this
 contradicts \eqref{contrconvdom}.

 The proof that $\|\Lambda_j^{\epsilon,1}\|_{L^2(\Omega\times D\times [0,T])}\to 0$
 is identical.
 The other terms $G ^{\epsilon,i}$,
 $\Lambda_j^{\epsilon,i}$ are treated in a standard way using
 Propositions
 \ref{lpfirstvariation} and \ref{lpsecondvariation}.
\qed

Define

$$
\begin{array}{lll}
\delta^\epsilon l(t,x)&=&
l(t,x,X_t(x),u^\epsilon_t)- l(t,x,X_t(x),u_t),
\\
\delta^\epsilon l'(t,x)&=&
l'(t,x,X_t(x),u^\epsilon_t)- l'(t,x,X_t(x),u_t).
\end{array}
$$

\begin{proposition}\label{costvariationprelim} We have
\begin{equation}\label{varcostprelim}
    \begin{array}{l}\dis
J(u^\epsilon)-J(u)=\E
\int_0^T\int_D \delta^\epsilon l(t,x)\, m(dx)\,dt
\\
+\dis
\E\int_0^T\int_D  l'(t,x,X_t(x),u_t)(
Y_t^\epsilon(x)+Z_t^\epsilon(x)) \,m(dx) \,dt
\\\dis
+\frac12
\E\int_0^T\int_D  l''(t,x,X_t(x),u_t)
Y_t^\epsilon(x)^2\,m(dx) \,dt
\\
+\dis
\E \int_D  h'( x,X_t(x) )(
Y_T^\epsilon(x)+Z_T^\epsilon(x)) \,m(dx)
+\frac12
 \E\int_D
 h''( x,X_T(x))Y_T^\epsilon(x)^2 \,m(dx)
+ o(\epsilon).
\end{array}
\end{equation}

\end{proposition}

\noindent{\bf Proof.} We still denote $R^\epsilon=Y^\epsilon+Z^\epsilon$.
We have
$$
\begin{array}{l}\dis
J(u^\epsilon)-J(u)=\E
\int_0^T\int_D [l(t,x,X_t^\epsilon(x),u^\epsilon_t)- l(t,x,X_t(x),u_t)]\, m(dx)\,dt
\\\dis\qquad
+\E\int_D
 [h( x,X_T^\epsilon(x))-h( x,X_T(x))]  \,m(dx).
\end{array}
$$
We first consider
$$
\E
\int_0^T\int_D [l(t,x,X_t^\epsilon(x),u^\epsilon_t)- l(t,x,X_t(x),u_t)]\, m(dx)\,dt
=A_1+A_2+A_3,
$$
where
\begin{eqnarray*}
  A_1 &=& \E
\int_0^T\int_D [l(t,x,X_t^\epsilon(x),u^\epsilon_t)
- l(t,x,X_t(x)+R_t^\epsilon(x),u_t^\epsilon)]\, m(dx)\,dt ,\\
  A_2 &=& \E
\int_0^T\int_D [
 l(t,x,X_t(x)+R_t^\epsilon(x),u_t^\epsilon)
 - l(t,x,X_t(x)+R_t^\epsilon(x),u_t)]\, m(dx)\,dt , \\
  A_3 &=& \E \int_0^T\int_D [
  l(t,x,X_t(x)+R_t^\epsilon(x),u_t)
  - l(t,x,X_t(x),u_t)]\, m(dx)\,dt .
\end{eqnarray*}
From Proposition \ref{variations} it follows that $A_1=o(\epsilon)$. Next applying
the Taylor formula \eqref{taylor} twice in $A_2$ we have
$$
\begin{array}{l}\dis
A_2= \E
\int_0^T\int_D \bigg( \delta^\epsilon l(t,x)+
\delta^\epsilon l'(t,x)R^\epsilon_t(x)
\\\dis+
\int_0^1\int_0^1 [ l''(t,x,X_t(x)+\lambda\mu R_t^\epsilon(x),u_t^\epsilon)
-  l''(t,x,X_t(x)+\lambda\mu R_t^\epsilon(x),u_t)]\,\mu d\lambda d\mu\,R^\epsilon_t(x) ^2
\bigg) \, m(dx)\,dt
\\\dis
= \E
\int_0^T\int_D   \delta^\epsilon l(t,x)\, m(dx)\,dt  +o(\epsilon),
\end{array}
$$
as it follows easily from Propositions \ref{lpsecondvariation} and \ref{lpsecondvariation}.
Applying
the Taylor formula \eqref{taylordue}  we have
$$
\begin{array}{l}\dis
A_3= \E
\int_0^T\int_D \bigg( l'(t,x,X_t(x),u_t)R_t^\epsilon(x)+
\frac12 l''(t,x,X_t(x),u_t)R_t^\epsilon(x)^2
\\\dis+
\int_0^1\int_0^1 [ l''(t,x,X_t(x)+\lambda\mu R_t^\epsilon(x),u_t)
-  l''(t,x,X_t(x),u_t)]\,\mu d\lambda d\mu\,R^\epsilon_t(x) ^2
\bigg) \, m(dx)\,dt
\\\dis
= \E
\int_0^T\int_D
\bigg( l'(t,x,X_t(x),u_t)R_t^\epsilon(x)+
\frac12 l''(t,x,X_t(x),u_t)Y_t^\epsilon(x)^2\bigg)
\, m(dx)\,dt  +o(\epsilon).
\end{array}
$$
The last equality is verified noting that
$$
\E
\int_0^T\int_D  l''(t,x,X_t(x),u_t)\,(2Y_t^\epsilon(x)Z_t^\epsilon(x)
+Z_t^\epsilon(x)^2)
\, m(dx)\,dt  =o(\epsilon),
$$
by Propositions \ref{lpsecondvariation} and \ref{lpsecondvariation}, and
that
$$
\E
\int_0^T\!\!\int_D \bigg(
\int_0^1\!\!\int_0^1 [ l''(t,x,X_t(x)+\lambda\mu R_t^\epsilon(x),u_t)
-  l''(t,x,X_t(x),u_t)]\,\mu d\lambda d\mu\,R^\epsilon_t(x) ^2
\bigg) \, m(dx)\,dt=o(\epsilon)
$$
which can be proved by the same arguments used to treat
the term $G^{\epsilon,1}$ in the proof of
Proposition \ref{variations}.

In a similar way one proves
$$
\begin{array}{l}\dis
\E\int_D
 [h( x,X_T^\epsilon(x))-h( x,X_T(x))]  \,m(dx)
 \\
 =\dis
\E \int_D  h'( x,X_t(x) )(
Y_T^\epsilon(x)+Z_T^\epsilon(x)) \,m(dx)
+\frac12
 \E\int_D
 h''( x,X_T(x))Y_T^\epsilon(x)^2 \,m(dx)
+ o(\epsilon),
\end{array}
$$
and the proof is finished.
\qed

\subsection{The first adjoint process}

The first adjoint process is defined as the solution of the backward stochastic
PDE
\begin{equation}\label{firstadjoint}
    \left\{\begin{array}{lll}
    -dp_t(x)&=&\dis -dq_t^j(x)\,dW^j_t
    +[
    A^*p_t(x)
     + b'(t,x,X_t(x),u_t)p_t(x)
     \\
    &&\dis
    +\sigma_j'(t,x,X_t(x),u_t)q_t^j(x)
    + l'(t,x,X_t(x),u_t)]\,dt
    \\
    p_T(x)&=&h'(x,X_T(x))
\end{array}\right.
\end{equation}
where $A^*$ denotes the adjoint of $A$ in $L^2$.
By the result in \cite{HuPe} there exists a unique solution,
i.e. a progressive
process $(p_t,q_t^1,\ldots, q_t^d )_{t\in [0,T]}$ with values in $(L^2)^{d+1}$, such that
$$
\sup_{t\in[0,T]}\E\|p_t\|_2^2+
\E\int_0^T\sum_{j=1}^d\|q^j_t\|_2^2\,dt
<\infty,
$$
and satisfying the equation
in the  mild sense:
 for every $t\in[0,T]$,
\begin{equation}\label{firstadjointmild}
\begin{array}{lll}\dis
    p_t+\int_t^Te^{(s-t)A^*}q_s^j \,dW^j_s &=&\dis
    e^{(T-t)A^*} h'(\cdot,X_T(\cdot))+
  \int_t^T e^{(s-t)A^*}[
    b'(s,\cdot,X_s(\cdot),u_s)p_s(\cdot)
    \\
&&\dis
    +
     \sigma_j'(s,\cdot,X_s(\cdot),u_s)q^j_s(\cdot)
    + l'(s,\cdot,X_s(\cdot),u_s)]
\,ds,
    \; \P-a.s.
    \end{array}
\end{equation}
where $(e^{tA^*})$ denotes the adjoint semigroup of $(e^{tA})$ in $L^2$,
which admits $A^*$ as its generator.

\begin{proposition}\label{costvariation}
Define
\begin{equation}\label{defHeh}
\begin{array}{lll}
    \bar H(t,x)&=& l''(t,x,X_t(x),u_t) +
p_t(x)  b''(t,x,X_t(x),u_t) + q^j_t(x)  \sigma_j''(t,x,X_t(x),u_t),
\\
\bar h(x)&=& h''( x,X_T(x)).
\end{array}
\end{equation}
Then we have
$$
\begin{array}{lll}\dis
J(u^\epsilon)-J(u)&=&\dis
\E
\int_0^T\int_D [
\delta^\epsilon l(t,x) +
p_t(x)\delta^\epsilon b(t,x) + q^j_t(x)\delta^\epsilon \sigma_j(t,x)
]\, ds\,m(dx)
\\
&&+\dis
\frac12
\E\int_0^T\int_D \bar H(t,x) \,Y_t^\epsilon(x)^2\, ds\,m(dx)
+
\frac12
 \E\int_D
 \bar h(x)\, Y_T^\epsilon(x)^2 \,m(dx)
+ o(\epsilon).
\end{array}
$$
\end{proposition}

\noindent{\bf Proof.}
We claim that the following duality relations hold:
\begin{equation}\label{dualityone}
\begin{array}{l}\dis
  \E\int_0^T\int_D  l'(t,x,X_t(x),u_t)
Y_t^\epsilon(x)\,m(dx) \,dt
+ \E \int_D  h'( x,X_t(x) )Y_T^\epsilon(x) \,m(dx)
\\\qquad
  \dis
=   \E\int_0^T\int_D  \delta^\epsilon\sigma_j(t,x)
q_t^j(x)\,m(dx) \,dt,
\end{array}
\end{equation}
\begin{equation}\label{dualitytwo}
\begin{array}{l}\dis
    \E\int_0^T\int_D  l'(t,x,X_t(x),u_t)
Z_t^\epsilon(x)\,m(dx) \,dt + \E \int_D  h'( x,X_t(x) )Z_T^\epsilon(x) \,m(dx)
\\\qquad
  \dis
=  \E\int_0^T\int_D  [\delta^\epsilon b(t,x)+
\frac12 b''(t,x,X_t(x),u_t) Y_t^\epsilon(x)^2+
\delta^\epsilon b'(t,x)
Y_t^\epsilon(x)]
p_t(x)\,m(dx) \,dt
\\\qquad
  \dis   +\E\int_0^T\int_D  [
\frac12 \sigma_j''(t,x,X_t(x),u_t) Y_t^\epsilon(x)^2+
\delta^\epsilon\sigma_j'(t,x)
Y_t^\epsilon(x)]
q_t^j(x)\,m(dx) \,dt.
\end{array}
\end{equation}
 If
$A$ is a bounded operator and equations \eqref{firstvariation} and
\eqref{secondvariation} are valid in the sense of Ito differentials in $L^2$
then \eqref{dualityone} and \eqref{dualitytwo}
 follow from an application of the Ito formula to the processes
$\<Y_t^\epsilon,p_t\>_{L^2}$ and $\<Z_t^\epsilon,p_t\>_{L^2}$ respectively,
where $\<\cdot,\cdot\>_{L^2}$ denotes the scalar product in $L^2$.
In the general case a regularization procedure is needed, where in particular
the operator $A$ is replaced by its Yosida approximation $A_n$ and then
$n\to\infty$. We omit writing down this standard part of the proof:
one can find the details of these   arguments (applied to BSDEs) in \cite{Te2}  or (applied to forward SDEs and control problems) in \cite{Te1}. One can also look at
Subsection \ref{technicalproof}
below where we use similar arguments in a more complicated setting.

 Now the proof is concluded substituting \eqref{dualityone} and \eqref{dualitytwo}
 in \eqref{varcostprelim}, provided we can prove
 $$
  \E\int_0^T\int_D
\delta^\epsilon b'(t,x)
Y_t^\epsilon(x)
p_t(x)\,m(dx) \,dt= o(\epsilon),
\quad
\E\int_0^T\int_D
\delta^\epsilon\sigma_j'(t,x)
Y_t^\epsilon(x)
q_t^j(x)\,m(dx) \,dt
= o(\epsilon).
$$
Since the proof is very similar, we
only prove the second equality. Since $\delta^\epsilon\sigma_j'$
is bounded and
 supported in $[t_0,t_0+\epsilon]$ we have, using the
 H\"older inequality
 and recalling the norm
$|||\cdot|||_p$ introduced in \eqref{normpprocess},
$$
 \begin{array}{l}\dis
 \left| \E\int_0^T\int_D
\delta^\epsilon\sigma_j'(t,x)
Y_t^\epsilon(x)
q_t^j(x)\,m(dx) \,dt
 \right|
 \le C\E\int_0^T
1_{[t_0,t_0+\epsilon]}(t)
\|Y_t^\epsilon\|_2
\|q_t\|_2  \,dt
\\\dis\qquad
 \le C |||Y^\epsilon|||_2
\int_0^T
1_{[t_0,t_0+\epsilon]}(t)
(\E\|q_t\|_2^2)^{1/2}  \,dt
\le C\, |||Y^\epsilon|||_2
\left(\int_0^T
1_{[t_0,t_0+\epsilon]}(t)
\E\|q_t\|_2^2 \,dt \right)^{1/2}\sqrt{\epsilon}.
 \end{array}
 $$
 The last integral tends to $0$ as $\epsilon\to 0$, since
 $\E\int_0^T \|q_t\|_2^2\,dt
<\infty$. It follows that
the right-hand side is $o(\epsilon)$, because
$|||Y^\epsilon|||_2\le C\,\sqrt{\epsilon}$
by
Proposition \ref{lpfirstvariation}.
\qed

\subsection{Some formal computations and heuristics}

In order to motivate some of the constructions below, and to make a connection
with the finite-dimensional case treated in \cite{Pe},
in this paragraph we proceed in a formal way.

We wish to prove that

\begin{equation}\label{finalduality}
    \begin{array}{l}\dis
\E\int_0^T\int_D \bar H(t,x) Y_t^\epsilon(x)^2\, ds\,m(dx)
+
 \E\int_D
 \bar h( x)Y_T^\epsilon(x)^2 \,m(dx)
 \\\dis  =
 \E\int_0^T
 \< P_t \delta^\epsilon \sigma_j(t,\cdot),\delta^\epsilon \sigma_j(t,\cdot)\>_{L^2}\,dt
+ o(\epsilon),
\end{array}
\end{equation}
for an appropriate operator-valued process $P_t$. In view of
Proposition \ref{costvariation} the stochastic maximum principle
then can be shown to hold by the usual arguments as in \cite{Pe} or \cite{YoZh}.

We denote by $H_t$ the multiplication
operator by the function $\bar H(t,\cdot)$ and
by $h$ the multiplication
operator by the function $\bar h(\cdot)$. We
pretend that they are bounded
 operators on the space $L^2$.

 Next we consider the operator-valued BSDE
\begin{equation}\label{BSDEoperatorvalued}
    \left\{\begin{array}{lll}
    -dP_t &=&\dis -Q_t^j \,dW^j_t
    +[
    A^*P_t + P_tA
     + B_tP_t +P_t B_t
    +C_t^j P_t C_t^j + C_t^jQ_t^j +  Q_t^jC_t^j
    + H_t ]\,dt
    \\
    P_T &=&h,
\end{array}\right.
\end{equation}
 where by $B_t$, $C^j_t$ we denote the (self-adjoint) multiplication operators
 by $b'(t,\cdot,X_t(\cdot),u_t)$ and
$ \sigma_j'(t,\cdot,X_t(\cdot),u_t)$ respectively.
Suppose that we can find a good solution in the space of bounded linear operators
on $L^2$. Then applying the Ito formula to $\<P_t Y^\epsilon_t,Y^\epsilon_t\>_{L^2}$,
integrating from $0$ to $T$  and taking expectations
we obtain
$$
\begin{array}{l}\dis
\E\int_0^T\<  H_t  Y_t^\epsilon ,Y_t^\epsilon \>_{L^2}\, dt
+
 \E\<
   h Y_T^\epsilon,  Y_T^\epsilon\>_{L^2}
=
 \E\int_0^T[
 \< P_t \delta^\epsilon \sigma_j(t,\cdot),\delta^\epsilon \sigma_j(t,\cdot)\>_{L^2}
  \\\dis
 +
2 \< P_t \delta^\epsilon b(t,\cdot),Y_t^\epsilon \>_{L^2}+
2 \< P_t  C^j_t Y^\epsilon_t ,\delta^\epsilon b(t,\cdot) \>_{L^2}
-
2 \< Q^j_t Y^\epsilon_t ,\delta^\epsilon \sigma_j(t,\cdot) \>_{L^2}
 ]\,dt.
\end{array}
$$
If we were able to prove, in analogy with the finite-dimensional case, that
$$
\E\int_0^T[
2 \< P_t \delta^\epsilon b(t,\cdot),Y_t^\epsilon \>_{L^2}+
2 \< P_t  C^j_t Y^\epsilon_t ,\delta^\epsilon b(t,\cdot) \>_{L^2}
-
2 \< Q^j_t Y^\epsilon_t ,\delta^\epsilon \sigma_j(t,\cdot) \>_{L^2}
 ]\,dt=
o(\epsilon),
$$
then \eqref{finalduality} would follow and the proof would be finished. However
in this argument finding a  solution
of the operator-valued BSDE
\eqref{BSDEoperatorvalued} that allows to make the previous
argument rigorous seems a very difficult task. So we follow a different
strategy of proof, that we outline below.

For fixed $t\in[0,T]$ and $f\in H$, denote by $(Y^{t,f}_s)_{s\in[t,T]}$ the mild solution to
$$
    \left\{\begin{array}{lll}
    dY^{t,f}_s(x)&=&\dis
    AY^{t,f}_s(x)\,ds + b'(s,x,X_s(x),u_s)Y^{t,f}_s(x)\,ds
     +  \sigma_j'(s,x,X_s(x),u_s)Y^{t,f}_s(x)\,dW^j_s,
    \\
    Y^{t,f}_t(x) &=&f(x)
\end{array}\right.
$$
This equation has to be compared with (\ref{firstvariation}).

Then taking $g\in {L^2}$,
applying the Ito formula to $\<P_s Y^{t,f}_s,Y^{t,g}_s\>_{L^2}$ over the interval
$[t,T]$,
integrating from $t$ to $T$  and taking conditional expectation given $\calf_t$
we formally obtain
\begin{equation}\label{bilinearone}
    \begin{array}{lll}\dis
 \< P_t f,g\>_{L^2}&=&\dis
\E^{\calf_t}\int_t^T\<  H_s  Y_s^{t,f} ,Y_s^{t,g} \>_{L^2}\, ds
+
 \E^{\calf_t}\<
   h Y_T^{t,f},  Y_T^{t,g}\>_{L^2}
   \\&=&\dis
\E^{\calf_t}\int_t^T\int_D \bar H(s,x)  Y_s^{t,f}(x) Y_s^{t,g}(x)\, m(dx)\, ds
+
 \E^{\calf_t}\int_D \bar
   h(x) Y_T^{t,f}(x)    Y_T^{t,g}(x) \, m(dx).
   \end{array}
\end{equation}
The interesting fact is that this formula can be used to
{\em define}  $P_t$: more precisely, in Proposition
 \ref{solauxiliaryy} below, we will prove that
if $f\in L^4$ then
 $(Y^{t,f}_s)_{s\in [t,T]}$ is
a progressive
process with values in $L^4$, satisfying
$$
\sup_{s\in[t,T]}(\E^{\calf_t}\|Y^{t,f}_s\|_4^4)^{1/4}\le C\|f\|_4.
$$
As a consequence we will show that
the right-hand side of (\ref{bilinearone})   defines a continuous bilinear form
on $L^4$   (or equivalently
 a linear bounded operator from $L^4$ to $L^{4/3}=(L^4)^*$) and we will set,
 for $f,g\in L^4$,
  $$
  \< P_t f,g\>= \E^{\calf_t}\int_t^T\int_D \bar H(s,x)  Y_s^{t,f}(x) Y_s^{t,g}(x)\, m(dx)\, ds
+
 \E^{\calf_t}\int_D \bar
   h(x) Y_T^{t,f}(x)    Y_T^{t,g}(x) \, m(dx).
 $$
Note that no reference to the BSDE (\ref{BSDEoperatorvalued}) is needed to give this definition.
Finally, it turns out that  (\ref{finalduality}) can be proved to hold  with this definition of $P_t$.

\section{End of the proof of the stochastic maximum principle}

As explained above, we are going to introduce
the second adjoint process, an appropriate operator-valued process $(P_t)$
that will allow us to conclude the proof of Theorem \ref{mainresult}.

Throughout this section we assume that Hypotheses
\ref{standing}, \ref{basisinl4} and
\ref{regularityinl4} are satisfied.

The  symbols
 $(e^{tA})_{t\ge 0}$ and $A$ will also denote
the restriction of the semigroup to the space $L^4$
and its
 infinitesimal generator in $L^4$. We need to recall some standard
 facts and constructions on analytic semigroups: see for instance
\cite{pa} or \cite{lu}.
 Without loss of generality we can assume
 that $A$ is boundedly invertible (if not, we replace $A$ by $A-cI$ for
 sufficiently large constant $c>0$ and we modify the drift coefficient
 accordingly). The domain of $A$ is endowed with the norm
 $\|f\|_{D(A)}:= \|Af\|_4$.
By the analyticity assumption, one can define the fractional powers $(-A)^\eta$
of $-A$ in a standard
way, for every $\eta\in (0,1)$. Each fractional power is a linear,
in general unbounded, operator in $L^4$, with domain denoted
$D(-A)^\eta$.
Endowed with the norm
$\|f\|_{D(-A)^\eta}:=\|(-A)^\eta f\|_4$, each space $D(-A)^\eta$ is a Banach space
and we have the continuous embeddings
$$D(A)\subset D(-A)^\eta \subset D(-A)^\rho \subset L^4,
\qquad 0<\rho<\eta <1.
$$
By analyticity, $e^{tA}(L^4)\subset D(A)$ for every $t> 0$,
and  for every $0<\eta <1$ there exist  constants $C_1,C_\eta>0$ such that
for every $f\in L^4$ and $t\in (0,T]$,
$$
\|e^{tA}f\|_{ D(A)}=\|Ae^{tA}f\|_4\le \frac{C_1}{t}\|f\|_4,
\qquad
\|e^{tA}f\|_{ D(-A)^\eta}=\|(-A)^\eta e^{tA}f\|_4\le \frac{C_\eta}{t^\eta}\|f\|_4.
$$
Finally, as a consequence of the compact embedding $D(A)\subset L^4$,
every embedding $D(-A)^\eta\subset L^4$ is also compact, $0<\eta <1$.

\begin{remark}{\em In most of what follows, we will only
use the estimate
$\|(-A)^\eta e^{tA}f\|_4\le \frac{C_\eta}{t^\eta}\|f\|_4$
and
the compact embedding $D(-A)^\eta\subset L^4$
for one, sufficiently small value of  $\eta >0$. This might
eventually lead to a weakening of Hypothesis \ref{regularityinl4},
but we not discuss those extensions in this paper.
}\end{remark}

For fixed $t\in[0,T]$ and $f\in L^4$, we consider the stochastic PDE
\begin{equation}\label{auxiliaryy}
    \left\{\begin{array}{lll}
    dY^{t,f}_s(x)&=&\dis
    AY^{t,f}_s(x)\,ds + b'(s,x,X_s(x),u_s)Y^{t,f}_s(x)\,ds
     +  \sigma_j'(s,x,X_s(x),u_s)Y^{t,f}_s(x)\,dW^j_s,
    \\
    Y^{t,f}_t(x) &=&f(x).
\end{array}\right.
\end{equation}
As a special case of Proposition \ref{lptrajectory} (with $p=4$), for every $t\in [0,T]$
there exists a unique mild solution, i.e. an adapted
process
$(Y^{t,f}_s)_{s\in[t,T]}$
with continuous
trajectories  in $L^4$, satisfying $\P$-a.s.
$$
    Y^{t,f}_s=
    e^{(s-t)A}f
         + \int_t^s e^{(s-r)A} b'(r,\cdot,X_r(\cdot),u_r)Y^{t,f}_r (\cdot)\,dr,
     + \int_t^s e^{(s-r)A} \sigma_j'(r,\cdot,X_r(\cdot),u_r)Y^{t,f}_r(\cdot) \,dW^j_r,
$$
for  every $s\in [t,T]$.
In addition we have
$
\sup_{0\le t\le s\le T} \E\|Y^{t,f}_s\|_4^4<\infty.
$

\begin{proposition}\label{solauxiliaryy}
There exists
a constant $C$ such that for $f\in L^4$, $0\le t\le s\le T$
 \begin{equation}\label{stimaunifauxiliaryy}
(\E^{\calf_t}\|Y^{t,f}_s\|_4^4)^{1/4}\le C\|f\|_4,\qquad \P-a.s.
 \end{equation}
and for  $0\le t\le t+h\le s\le T$
 \begin{equation}\label{stimaunifauxiliaryydue}
(\E\|Y^{t+h,f}_s-Y^{t,f}_s\|_4^4)^{1/4}\le C[
\sup_{t\in [0,T]}\|(e^{tA}- e^{(t+h)A})f\|_4+h^{1/2}\|f\|_4].
 \end{equation}

Moreover for every $\eta\in (0,1/4)$ there exists
a constant $C_\eta $ such that for $f\in D(-A)^\eta\subset L^4$, $0\le t< s\le T$
 \begin{equation}\label{stimaunifauxiliaryytre}
(\E^{\calf_t}\|Y^{t,(-A)^\eta f}_s\|_4^4)^{1/4}\le C_\eta (s-t)^{-\eta}\|f\|_4,\qquad \P-a.s.
 \end{equation}
We notice that the above relation indicates that equation (\ref{auxiliaryy}) regularizes the initial data
(roughly speaking sends data in $ D(-A)^{-\eta}$ to $L^4$.
\end{proposition}
{\bf Proof.}
For brevity we write the proof in the case  $b\equiv 0$ and denote by
$ C_j(r)$ the multiplication operator in $L^4$ by the
(bounded) function
$\sigma_j'(r,\cdot,X_r(\cdot),u_r)$. So the equation for $Y^{t,x}$
is
$$
    Y^{t,f}_s=
    e^{(s-t)A}f
     + \int_t^s e^{(s-r)A} C_j(r)Y^{t,f}_r \,dW^j_r,
    \qquad s\in [t,T].
$$
Using the conditional inequality (\ref{isometryiitotercond}) for $p=4$ we obtain
$$
  \E^{\calf_t}\left\| \int_t^s e^{(s-r)A} C_j(r)Y^{t,f}_r \,dW^j_r\right\|_4^4\le
  C
   \int_t^s \E^{\calf_t}\|e^{(s-r)A} C_j(r)Y^{t,f}_r\|_4^4 \,dr
   \le C
   \int_t^s \E^{\calf_t}\|Y^{t,f}_r\|_4^4 \,dr,
$$
and since $\|e^{(s-t)A}f\|_4^4\le C\|f\|_4^4$ it follows that
for every $s\in [t,T]$ we have, $\P$-a.s.
\begin{equation}\label{stimecondizionali}
    \E^{\calf_t}\|Y^{t,f}_s\|_4^4\le
    C\|f\|_4^4
     +C
   \int_t^s \E^{\calf_t}\|Y^{t,f}_r\|_4^4 \,dr.
\end{equation}
Take a dense countable set  $D\subset [t,T]$. Then, $\P$-a.s.,
\eqref{stimecondizionali} holds simultaneously  for every $s\in D$.
Since $Y^{t,x}$ has continuous trajectories in $L^4$, there exists
a set $N$ with  $\P(N)=0$  such that
$s\to \|Y^{t,f}_s(\omega)\|_4^4$ is continuous on $[t,T]$ for
every $\omega\notin N$.
Discarding a set of $\P$-measure zero, and
given any $s\in [t,T]$, we take a sequence $(s_n)\subset D$,
$s_n\to s$ and by the conditional Fatou Lemma
$$
\begin{array}{l}\dis
  \E^{\calf_t}\|Y^{t,f}_s\|_4^4=
  \E^{\calf_t}\liminf_{n\to\infty}\|Y^{t,f}_{s_n}\|_4^4
  \le\liminf_{n\to\infty}
  \E^{\calf_t}\|Y^{t,f}_{s_n}\|_4^4
  \\\dis
  \le
    C\|f\|_4^4
     +C\liminf_{n\to\infty}
   \int_t^{s_n} \E^{\calf_t}\|Y^{t,f}_r\|_4^4 \,dr=
    C\|f\|_4^4
     +C
   \int_t^{s} \E^{\calf_t}\|Y^{t,f}_r\|_4^4 \,dr.
   \end{array}
   $$

It follows  that, $\P$-a.s. \eqref{stimecondizionali} holds
for every $s\in [t,T]$, so that
(\ref{stimaunifauxiliaryy}) follows from
a pathwise application of Gronwall's lemma.

 The proof of (\ref{stimaunifauxiliaryytre}) is very similar: by
 (\ref{stimecondizionali}) we have
$$
\begin{array}{l}\dis
    \E^{\calf_t}\|Y^{t,(-A)^\eta f}_s\|_4^4\le
    C\|(-A)^\eta e^{(s-t)A}f\|_4^4
     +C
   \int_t^s \E^{\calf_t}\|Y^{t,(-A)^\eta f}_r\|_4^4 \,dr
   \\\dis
   \le
    C_\eta (s-t)^{-4\eta}\|f\|_4^4
     +C
   \int_t^s \E^{\calf_t}\|Y^{t,(-A)^\eta f}_r\|_4^4 \,dr
\end{array}
$$
and (\ref{stimaunifauxiliaryytre}) follows again from
a variant of Gronwall's lemma.

To prove (\ref{stimaunifauxiliaryydue}) we first write, for $s\in [t+h,T]$,
$$
\begin{array}{l}\dis
     Y^{t+h,  f}_s-Y^{t,  f}_s= (e^{(s-t-h)A}- e^{(s-t)A})f -
   \int_t^{t+h} e^{(s-r)A}C_j(r)  Y^{t,  f}_r \,dW^j_r
   \\\dis
   +
    \int_{t+h}^s e^{(s-r)A}C_j(r)  (Y^{t+h,  f}_r -Y^{t,  f}_r )\,dW^j_r
     =:I+II+III.
\end{array}
$$
Then we have, using (\ref{isometryiitoter}) for $p=4$,
$$
\|I\|_4\le \sup_{t\in [0,T]}\|(e^{tA}- e^{(t+h)A})f\|_4,
$$
$$
\E\|II\|_4^4\le  c h\int_t^{t+h}\E \| Y^{t,  f}_r\|_4^4 \,dr
\le c h^2\| f\|_4^4,
$$
$$
\E\|III\|_4^4\le  c \int_{t+h}^s\E \| Y^{t+h,  f}_r-Y^{t,  f}_r\|_4^4 \,dr.
$$
Therefore
$$
\E\|     Y^{t+h,  f}_s-Y^{t,  f}_s\|_4^4\le
c[\sup_{t\in [0,T]}\|(e^{tA}- e^{(t+h)A})f\|_4 + h^2\| f\|_4^4]
+ c \int_{t+h}^s\E \| Y^{t+h,  f}_r-Y^{t,  f}_r\|_4^4 \,dr
$$
and (\ref{stimaunifauxiliaryydue}) follows  from Gronwall's lemma.
 \qed

Recall that we denoted by $\call$
 the space of linear
bounded operators $L^4\to (L^4)^*=L^{4/3}$ endowed with the usual operator norm
and with the Borel $\sigma$-algebra of the weak topology.
The duality between $g\in L^4$ and $h\in L^{4/3}$
is denoted $\<h,g\>$.
We   note that, by the H\"older inequality,
 every $H\in L^2$ can be identified with the corresponding multiplication
 operator, i.e. with a unique $H\in \call$ satisfying
    $$\< H f,g\>= \int_D H(x)  f(x)g(x)\, m(dx), \qquad
    f,g\in L^4$$ and,
    moreover, $\|H\|_\call \le \|H\|_2$.
Similar remarks apply to
    $\bar H(t,x)$ and $
\bar h(x)$  defined in
(\ref{defHeh}).

The definition of the second adjoint process $P$, along with some of
its properties, is given in the following proposition.

\begin{proposition}\label{processp}
There exists a progressive process $(P_t)_{t\in [0,T]}$ with values
in $\call$, such that for $t\in [0,T]$, $f,g\in L^4$,
  \begin{equation}\label{defpt}
    \< P_t f,g\>= \E^{\calf_t}\int_t^T\<\bar H_s Y_s^{t,f}, Y_s^{t,g}\>\, ds
+
 \E^{\calf_t}\< \bar
   h  Y_T^{t,f},    Y_T^{t,g}  \>, \qquad \P-a.s.
  \end{equation}

We have
  \begin{equation}\label{stimainelleperp}
    \sup_{t\in [0,T]}\E\|P_t\|_\call^2<\infty,
  \end{equation}
and for every $f,g\in L^4$ we have, for $\epsilon \downarrow 0$,
  \begin{equation}\label{continelleperp}
\E|\<P_{t+\epsilon}-P_t)f,g\>|\to 0.
  \end{equation}

  Moreover,
  for every $\eta\in (0,1/4)$ there exists
a constant $C_\eta $ such that for $f,g\in D(-A)^\eta\subset L^4$, $0\le t<  T$,
 \begin{equation}\label{stimaindaetaperp}
| \<P_t(-A)^\eta f,(-A)^\eta g\>| \le
C_\eta \|f\|_4\|g\|_4 (T-t)^{-2\eta} \left[
\left(\int_t^T\E^{\calf_t}\|\bar H_s\|_2^2ds\right)^{1/2}
+
\left(\E^{\calf_t}\|\bar h\|_2^2 \right)^{1/2}
\right]
,\; \P-a.s.
 \end{equation}
which immediately implies
 \begin{equation}\label{stimaindaetaperpdue}
 \begin{array}{l}
\E \sup\Big\{
| \<P_t(-A)^\eta f,(-A)^\eta g\>|^2\;:\;
f,g\in D(-A)^\eta,
\|f\|_4\le 1,\|g\|_4\le 1
\Big\}
\\\dis
 \le
C_\eta  (T-t)^{-4\eta} \left[\E
\int_0^T \|\bar H_s\|_2^2ds
+
 \E \|\bar h\|_2^2 \right].
 \end{array}
 \end{equation}

   \end{proposition}

\begin{remark}\begin{em}
\begin{enumerate}
\item
Formula (\ref{defpt}) can be written more explicitly as follows:
$$
    \< P_t f,g\>= \E^{\calf_t}\int_t^T\int_D \bar H(s,x)  Y_s^{t,f}(x) Y_s^{t,g}(x)\, m(dx)\, ds
+
 \E^{\calf_t}\int_D \bar
   h(x) Y_T^{t,f}(x)    Y_T^{t,g}(x) \, m(dx).
$$
Clearly, (\ref{defpt})  defines uniquely $(P_t)$ up to modification.
\item
In the following, for $T\in\call$ we will use the notation
\begin{equation}\label{normextension}
    |||T||| :=
\sup\Big\{
| \<T(-A)^\eta f,(-A)^\eta g\>|\;:\;
f,g\in D(-A)^\eta,
\|f\|_4\le 1,\|g\|_4\le 1
\Big\}
\end{equation}
(\ref{stimaindaetaperpdue}) can then be written
\begin{equation}\label{stimaindaetaperptre}
\E\, ||| P_t|||^2
 \le
C_\eta  (T-t)^{-4\eta} \left[\E
\int_0^T \|\bar H_s\|_2^2ds
+
 \E \|\bar h\|_2^2 \right].
 \end{equation}
 \item Note that for $L^4$-valued, $\calf_t$-measurable random variables $F,G$ we have
 \begin{equation}\label{ptcomposition}
    \< P_t F,G\>= \E^{\calf_t}\int_t^T\<\bar H_s Y_s^{t,F}, Y_s^{t,G}\>\, ds
+
 \E^{\calf_t}\< \bar
   h  Y_T^{t,F},    Y_T^{t,G}  \>, \qquad \P-a.s.
  \end{equation}
  The equality being trivial if $F$ and $G$ are simple random variables and easily passing to the limit.
\end{enumerate}
\end{em}
\end{remark}

{\bf Proof of Proposition \ref{processp}.}
Fix
  $\eta\in (0,1/4)$, $f,g\in D(-A)^\eta\subset L^4$, $0\le t<  T$.
Using the conditional H\"older inequality
and (\ref{stimaunifauxiliaryytre}) we have
\begin{equation}\label{stimaindaetaperpdueprelim}
 \begin{array}{l}\dis
\left|  \E^{\calf_t}\int_t^T\<\bar H_s Y_s^{t,(-A)^\eta f}, Y_s^{t,(-A)^\eta g}\>\, ds
+
 \E^{\calf_t}\< \bar
   h  Y_T^{t,(-A)^\eta f},    Y_T^{t,(-A)^\eta g}  \>\right|
   \\\dis
\le
 \E^{\calf_t}\int_t^T\|\bar H_s\|_2\| Y_s^{t,(-A)^\eta f}\|_4\| Y_s^{t,(-A)^\eta g}\|_4\, ds
+
 \E^{\calf_t}[\| \bar
   h\|_2\|  Y_T^{t,(-A)^\eta f}\|_4\|    Y_T^{t,(-A)^\eta g}  \|_4]
\\\dis
\le
 \int_t^T(\E^{\calf_t}\|\bar H_s\|_2^2)^{1/2}
 (\E^{\calf_t}\| Y_s^{t,(-A)^\eta f}\|_4^4)^{1/4}
 (\E^{\calf_t}\| Y_s^{t,(-A)^\eta g}\|_4^4)^{1/4}\, ds
\\\dis \qquad
+
( \E^{\calf_t}\| \bar
   h\|_2\|^2)^{1/2}
   (\E^{\calf_t}\| Y_T^{t,(-A)^\eta f}\|_4^4)^{1/4}
   (\E^{\calf_t}\|    Y_T^{t,(-A)^\eta g}  \|_4^4)^{1/4}
   \\\dis
\le c\|   f \|_4 \|   g \|_4
 \int_t^T(\E^{\calf_t}\|\bar H_s\|_2^2)^{1/2}
 (s-t)^{-2\eta}\, ds
+c\|   f \|_4 \|   g \|_4
( \E^{\calf_t}\| \bar    h\|_2\|^2)^{1/2}
(T-t)^{-2\eta}
\\\dis
\le c\|   f \|_4 \|   g \|_4\left[
\left( \int_t^T \E^{\calf_t}\|\bar H_s\|_2^2ds\right)^{1/2}
\left( \int_t^T  (s-t)^{-4\eta} ds\right)^{1/2}
+( \E^{\calf_t}\| \bar    h\|_2\|^2)^{1/2}
(T-t)^{-2\eta}\right]
\\\dis
 \le c
 \|f\|_4\|g\|_4 (T-t)^{-2\eta} \left[
\left(\int_t^T\E^{\calf_t}\|\bar H_s\|_2^2ds\right)^{1/2}
+
\left(\E^{\calf_t}\|\bar h\|_2^2 \right)^{1/2}
\right],
 \end{array}
 \end{equation}
 where $c$ is a constant independent of $f,g,t$.
 Using (\ref{stimaunifauxiliaryy}) instead of
 (\ref{stimaunifauxiliaryytre})  this inequality also holds for
   $\eta=0$.

Now fix a dense set $F$ in $L^4$. For $f,g\in F$
let us define $\<P_tf,g\>$ by
formula (\ref{defpt}), by fixing an arbitrary version
of the conditional expectations on the right-hand side.
By (\ref{stimaindaetaperpdueprelim}) with $\eta=0$, there exists
a set $N$ of probability zero such that for $\omega\notin F$
we have
$$|\<P_t(\omega) f,g\>|\le c
 \|f\|_4\|g\|_4, \qquad f,g\in F. $$
 Thus, the mapping
 $(f,g)\mapsto  \<P_t(\omega)f,g\>$ extends from $F\times F$ to a continuous
 bilinear form on $L^4$ (or equivalently
 an element of $\call$), still denoted $P_t(\omega)$. Set
 $P_t(\omega)=0$ for $\omega\in N$. Using again
(\ref{stimaindaetaperpdueprelim}) with $\eta=0$, it is easily proved that
equality (\ref{defpt}) holds for every $f,g\in L^4$,
by approximating $f,g$ with elements of $F$.
Thus, an $\call$-valued process
   $(P_t)_{t\in [0,T]}$ has been constructed with the required properties.
   $(P_t)$ is adapted by construction.
Similar arguments also show the existence of a progressive
modification of $(P_t)$, as  required.

(\ref{stimaindaetaperp}) follows at once from (\ref{stimaindaetaperpdueprelim}).
(\ref{stimaindaetaperpdueprelim}) with $\eta=0$ gives
$$
\|P_t\|\le
c \,\left[
\left(\int_t^T\E^{\calf_t}\|\bar H_s\|_2^2ds\right)^{1/2}
+
\left(\E^{\calf_t}\|\bar h\|_2^2 \right)^{1/2}
\right],
$$
which implies  (\ref{stimainelleperp}).

It remains to prove (\ref{continelleperp}). We sketch the proof in the case
$\bar h=0$ for short.
$$\begin{array}{l}\dis
\<(P_{t+\epsilon}-P_t)f,g\>=
(\E^{\calf_{t+\epsilon}}- \E^{\calf_t})\int_t^T\<\bar H_s Y_s^{t,f}, Y_s^{t,g}\>\, ds
\\\dis
-\E^{\calf_{t+\epsilon}}\int_t^{t+\epsilon}\<\bar H_s Y_s^{t,f}, Y_s^{t,g}\>\, ds
+\E^{\calf_t}\int_{t+\epsilon}^T[\<\bar H_s Y_s^{{t+\epsilon},f}, Y_s^{{t+\epsilon},g}\>-
\<\bar H_s Y_s^{t,f}, Y_s^{t,g}\>]\, ds.
\end{array}
$$
The first summand tends to zero in $L^1(\Omega,\P)$ by the downwards martingale
convergence theorem, the third one due to (\ref{stimaunifauxiliaryydue}) and the second one
is easy to treat by dominated convergence Theorem.
\qed

We are now ready to finish the proof of our main result, by showing that
the formula (\ref{finalduality}) introduced during our heuristic
discussion actually holds (more precisely we will prove \eqref{finaldualitydue} below).

\noindent {\bf End of the proof of Theorem \ref{mainresult}.}
We claim that the following holds:
\begin{equation}\label{finaldualitydue}
\E\int_0^T\<\bar H_s  Y_s^\epsilon,Y_s^\epsilon \> \,ds
+
 \E\<
 \bar h Y_T^\epsilon,Y_T^\epsilon \>
 =
 \E\int_0^T
 \< P_s \delta^\epsilon \sigma_j(s,\cdot),\delta^\epsilon \sigma_j(s,\cdot)\>\,ds
+ o(\epsilon).
\end{equation}

Admitting this for a moment, if follows from
Proposition \ref{costvariation} that
$$
\begin{array}{lll}\dis
J(u^\epsilon)-J(u)&=&\dis
\E
\int_0^T\int_D [
\delta^\epsilon l(t,x) +
p_t(x)\delta^\epsilon b(t,x) + q^j_t(x)\delta^\epsilon \sigma_j(t,x)
]\, ds\,m(dx)
\\
&&\dis+
 \E\int_0^T
 \< P_s \delta^\epsilon \sigma_j(s,\cdot),\delta^\epsilon \sigma_j(s,\cdot)\>\,ds
+ o(\epsilon).
\end{array}
$$
The optimality of $u$ implies that $J(u^\epsilon)-J(u)\ge 0$. Diving by $\epsilon$
and letting $\epsilon\to 0$, the required conclusion is obtained by standard arguments,
see e.g. \cite{Pe} or \cite{YoZh}.

So it only remains to prove \eqref{finaldualitydue}.
Recalling that $Y_s^\epsilon=0$ for $s\le t_0$, the left-hand side of
(\ref{finaldualitydue}) equals
$$
\E\int_{t_0}^{t_0+\epsilon}\<\bar H_s  Y_s^\epsilon,Y_s^\epsilon \> \,ds
+
\E\int_{t_0+\epsilon}^T\<\bar H_s  Y_s^\epsilon,Y_s^\epsilon \> \,ds
+
 \E\<
 \bar h Y_T^\epsilon,Y_T^\epsilon \>.
$$
It is easily checked that the first integral is $o(\epsilon)$. Using the formula
$$
Y_s^\epsilon=
Y_s^{t_0+\epsilon,Y_{t_0+\epsilon}^\epsilon}, \qquad s\ge {t_0+\epsilon},
$$
which follows by comparing the equations
(\ref{auxiliaryy}) and
(\ref{firstvariation})
satisfied by $Y^{t_0+\epsilon,f}$
and $Y^\epsilon$,
we obtain
\begin{equation}\label{finaldualitytre}
    \begin{array}{l}\dis
\E\int_0^T\<\bar H_s  Y_s^\epsilon,Y_s^\epsilon \> \,ds
+
 \E\<
 \bar h Y_T^\epsilon,Y_T^\epsilon \>
 \\\dis
 = o(\epsilon)+
 \E\int_{t_0+\epsilon}^T\<\bar H_s
 Y_s^{t_0+\epsilon,Y_{t_0+\epsilon}^\epsilon}
,Y_s^{t_0+\epsilon,Y_{t_0+\epsilon}^\epsilon} \> \,ds
+
 \E\<
 \bar h Y_T^{t_0+\epsilon,Y_{t_0+\epsilon}^\epsilon}
 ,Y_T^{t_0+\epsilon,Y_{t_0+\epsilon}^\epsilon}\>
 \\\dis  = o(\epsilon)+
 \E
 \< P_{t_0+\epsilon}Y_{t_0+\epsilon}^\epsilon,
 Y_{t_0+\epsilon}^\epsilon\>,
\end{array}
\end{equation}
where the last equality follows from an application of
(\ref{ptcomposition}).
Next we claim that
\begin{equation}\label{deltapyy}
  \E  \< (P_{t_0+\epsilon}-P_{t_0})\, Y_{t_0+\epsilon}^\epsilon,
Y_{t_0+\epsilon}^\epsilon\>= o(\epsilon),
\end{equation}
\begin{equation}\label{deltapyydue}
  \E  \<  P_{t_0}  Y_{t_0+\epsilon}^\epsilon,
Y_{t_0+\epsilon}^\epsilon\>= \E\int_{t_0}^{t_0+\epsilon}
 \< P_s \delta^\epsilon \sigma_j(s,\cdot),\delta^\epsilon \sigma_j(s,\cdot)\>\,ds
+ o(\epsilon).
\end{equation}
The required formula
(\ref{finaldualitydue}) will now be a consequence of \eqref{deltapyy}
and \eqref{deltapyydue}, which are proved in the following
two subsection below. The proof of Theorem \ref{mainresult}
will then be finished.
\qed

\subsection{Proof of (\ref{deltapyy})}

It is convenient to rewrite (\ref{deltapyy}) in the form
\begin{equation}\label{deltapyybis}
  \E  \< (P_{t_0+\epsilon}-P_{t_0})\, \epsilon^{-1/2}Y_{t_0+\epsilon}^\epsilon,
\epsilon^{-1/2} Y_{t_0+\epsilon}^\epsilon\>\to 0. \end{equation}
By Proposition  \ref{lpfirstvariation} there exists a constant  $C_0$
independent of $\epsilon$ such that
\begin{equation}\label{stimeunifinepsperyeps}
    (\E\|\epsilon^{-1/2}Y^\epsilon_{t_0+\epsilon}\|_4^4)^{1/4}\le C_0,
\quad
(\E\|\epsilon^{-1/2}Y^\epsilon_{t_0+\epsilon}\|_4^8)^{1/8}\le C_0.
\end{equation}
Next we fix $\eta\in (0,1/4)$ and notice that for
every $\delta>0$ we have, by the Markov inequality,
$$
\P(
\|\epsilon^{-1/2}(-A)^{-\eta} Y^\epsilon_{t_0+\epsilon}\|_{D(-A)^\eta}>C_0\delta^{-1/4})
=
\P(
\|\epsilon^{-1/2}Y^\epsilon_{t_0+\epsilon}\|_4>C_0\delta^{-1/4})\le \delta.
$$
Therefore setting
$K_\delta =\{ f\in L^4\;:\; f\in D(-A)^\eta,
\|f\|_{D(-A)^\eta}\le C_0\delta^{-1/4}\}$
and denoting $\Omega_{\delta,\epsilon}$ the event
$\{\epsilon^{-1/2}(-A)^{-\eta} Y^\epsilon_{t_0+\epsilon}\in K_\delta\}$
we obtain
$$
\P(\Omega_{\delta,\epsilon}^c)=
\P(
\epsilon^{-1/2}(-A)^{-\eta} Y^\epsilon_{t_0+\epsilon}\notin K_\delta)
\le \delta.
$$
We note that, since $D(-A)^\eta$ is compactly embedded in $L^4$, the set $K_\delta$
is a compact subset of $L^4$. Moreover, for $f\in K_\delta $ we have
\begin{equation}\label{radiuskdelta}
    \|f\|_4\le c \|f\|_{D(-A)^\eta} \le cC_0 \delta^{-1/4},
\end{equation}
i.e. $K_\delta$ is contained in a ball of $L^4$ centered at $0$ with radius proportional to
$\delta^{-1/4}$.

We have
$$\begin{array}{l}
 \E  \< (P_{t_0+\epsilon}-P_{t_0})\, \epsilon^{-1/2}Y_{t_0+\epsilon}^\epsilon,
\epsilon^{-1/2} Y_{t_0+\epsilon}^\epsilon\>
\\\dis =
 \E  [\< (P_{t_0+\epsilon}-P_{t_0})\, \epsilon^{-1/2}Y_{t_0+\epsilon}^\epsilon,
\epsilon^{-1/2} Y_{t_0+\epsilon}^\epsilon\>1_{\Omega_{\delta,\epsilon}^c}]+
 \E  [\< (P_{t_0+\epsilon}-P_{t_0})\, \epsilon^{-1/2}Y_{t_0+\epsilon}^\epsilon,
\epsilon^{-1/2} Y_{t_0+\epsilon}^\epsilon\>1_{\Omega_{\delta,\epsilon}}]
\\\dis
=:A_1^\epsilon+ A_2^\epsilon.
\end{array}
$$
By the H\"older inequality
$$
|A_1^\epsilon|\le (\E\|P_{t_0+\epsilon}-P_{t_0}\|^2_\call)^{1/2}
(\E\|\epsilon^{-1/2}Y^\epsilon_{t_0+\epsilon}\|_4^8)^{1/4}\P(\Omega_{\delta,\epsilon}^c)^{1/4},
$$
and from (\ref{stimainelleperp}), (\ref{stimeunifinepsperyeps})
 we conclude that
 $
|A_1^\epsilon|\le c\P(\Omega_{\delta,\epsilon}^c)^{1/4} \le c \delta^{1/4}$
for some constant $c$ independent of $\delta$ and $\epsilon$.

On the other hand, recalling the definition of $\Omega_{\delta,\epsilon}$,
$$
|A_2^\epsilon|\le
\E  \sup_{f\in K_\delta }|\< (P_{t_0+\epsilon}-P_{t_0})\, (-A)^\eta f,
(-A)^\eta f\>1_{\Omega_{\delta,\epsilon}}|.
$$
Since $K_\delta$ is compact in $L^4$, it can be covered by a finite number $N_\delta$
of open balls with radius $\delta$ and centers denoted $f_i^\delta$, $i=1,\ldots,N_\delta$.
Since $D(-A)^\eta$ is dense in $L^4$, we can assume that $f_i^\delta\in D(-A)^\eta$.
Given $f\in K_\delta$, let $i$ be such that $\|f-f_i^\delta\|_4<\delta$; then writing
$$\begin{array}{l}
  \< (P_{t_0+\epsilon}-P_{t_0}) (-A)^\eta f,(-A)^\eta f\>=
  \< (P_{t_0+\epsilon}-P_{t_0}) (-A)^\eta f_i^\delta,(-A)^\eta f_i^\delta\>
\\\dis -
  \< (P_{t_0+\epsilon}-P_{t_0}) (-A)^\eta (f-f_i^\delta),(-A)^\eta (f-f_i^\delta)\>
  +2 \< (P_{t_0+\epsilon}-P_{t_0}) (-A)^\eta f,(-A)^\eta (f-f_i^\delta)\>
\end{array}
$$
and recalling the notation introduced in (\ref{normextension}) we obtain
$$\begin{array}{l}
|  \< (P_{t_0+\epsilon}-P_{t_0}) (-A)^\eta f,(-A)^\eta f\>|\le
  |\< (P_{t_0+\epsilon}-P_{t_0}) (-A)^\eta f_i^\delta,(-A)^\eta f_i^\delta\>|
\\\dis +
  |||P_{t_0+\epsilon}-P_{t_0}|||\,\delta^2
  +2 |||P_{t_0+\epsilon}-P_{t_0}|||\,\|f\|_4\,\delta.
\end{array}
$$
Recalling (\ref{radiuskdelta}) we conclude that
$$\begin{array}{l}\dis
\sup_{f\in K_\delta }|\< (P_{t_0+\epsilon}-P_{t_0})\, (-A)^\eta f,
(-A)^\eta f\>\le
\sum_{i=1}^{N_\delta}
 |\< (P_{t_0+\epsilon}-P_{t_0}) (-A)^\eta f_i^\delta,(-A)^\eta f_i^\delta\>|
\\\dis +
2\sup_{t\in[t_0,t_0+\epsilon]}|||P_{t}|||\,\delta^2
  +c\sup_{t\in[t_0,t_0+\epsilon]}||| P_{t }|||\,\delta^{3/4},
\end{array}
$$
for some constant $c$.
Taking expectation,
it follows from (\ref{stimaindaetaperptre})
that
$$
|A_2^\epsilon|\le
\sum_{i=1}^{N_\delta}
\E |\< (P_{t_0+\epsilon}-P_{t_0}) (-A)^\eta f_i^\delta,(-A)^\eta f_i^\delta\>|
+c(T-t_0-\epsilon)^{-2\eta}
[\delta^2
  + \delta^{3/4}],
$$
for some constant $c$ independent of $\epsilon$ and $\delta$. By
(\ref{continelleperp}) we conclude that
$$
\limsup_{\epsilon \downarrow 0} |A_2^\epsilon|\le
 c(T-t_0)^{-2\eta}
[\delta^2
  + \delta^{3/4}].
$$
Letting $\delta\to 0$ we obtain $|A_1^\epsilon|+|A_2^\epsilon|\to 0$ and the proof
of (\ref{deltapyy})
is finished.

\subsection{Proof of (\ref{deltapyydue})}
\label{technicalproof}

In order to make appropriate computations on
$\E  \< P_{t_0}\, Y_{t_0+\epsilon}^\epsilon,
Y_{t_0+\epsilon}^\epsilon\>$
we perform an approximation of  both $ P_{t_0}$ and $Y_{t_0+\epsilon}^\epsilon$.

To approximate $ P_{t_0}$
  we use the basis $(e_i)_{i\ge 1}$  of Hypothesis \ref{basisinl4}.
We introduce the projection operators $\Pi_Nf=\sum_{i=1}^N\<f,e_i\>_2e_i$, $f\in L^2$,
where $\<\cdot,\cdot\>_2$ denotes the scalar product of $L^2$. Each $\Pi_N$
is an orthogonal projection in $L^2$. Since we assume that $(e_i)_{i\ge 1}$
is a Schauder basis of $L^4$, the restriction
of $\Pi_N$ to $L^4$ is a bounded linear operator in $L^4$, satisfying
$\|\Pi_Nf-f\|_4\to 0$ for every $f\in L^4$  and $\sup_N\|\Pi_N\|_{L(L^4,L^4)}<\infty$.
Then we define
$$
P_t^N(\omega)f:
=\sum_{i,j=1}^N \<P_t(\omega)e_i,e_j\> \<e_i,f\>_2 e_j,
\qquad f\in L^4.
$$
Then $P^N_t(\omega)$ is a  linear bounded operator on $L^4$, which extends to a linear bounded operator
on $L^2$, with values in the finite-dimensional subspace spanned by $e_1,\ldots,e_N$.
Moreover
\begin{equation}\label{cambionorme}
    \<P_t^N(\omega)f,g\>_2
=\sum_{i,j=1}^N \<P_t(\omega)e_i,e_j\> \<e_i,f\>_2\<e_j,g\>_2
=\<P_t(\omega)\Pi_Nf,\Pi_Ng\>,
\qquad f,g\in L^4.
\end{equation}
In the following we will consider $P^N$ as a stochastic process with values in $\call_2(L^2)$,
the space of Hilbert-Schmidt operators on $L^2$.

In order to approximate $Y_{t_0+\epsilon}^\epsilon$ we introduce
$$
J_n=(nI-A)^{-1},\qquad A_n=AJ_n,\qquad Y_t^{\epsilon,n}(\omega)=J_nY_t^{\epsilon}(\omega).
$$
Note that $A_n$ are the Yosida approximations of the operator $A$.

We are going to approximate $\E  \< P_{t_0}\, Y_{t_0+\epsilon}^\epsilon,
Y_{t_0+\epsilon}^\epsilon\>$
by $\E  \< P_{t_0}^N\, Y_{t_0+\epsilon}^{\epsilon,n},
Y_{t_0+\epsilon}^{\epsilon,n}\>_2$.

$Y^{\epsilon,n}$ is a process with values in $L^2$ which admits
an Ito differential that we are going to compute. Recall
equation (\ref{firstvariationmild}) satisfied by $Y^\epsilon$, that we
now re-write in the following way: for $s\ge t_0$,
$$
    Y^\epsilon_s=
  \int_{t_0}^s e^{(s-r)A}[
    B(r)Y^\epsilon_r
    + \delta^\epsilon b(r)]\,dr
    + \int_{t_0}^s e^{(s-r)A}
    [ C_j(r)Y^\epsilon_r
    + \delta^\epsilon \sigma_j (r)]
\,dW^j_r,
    \qquad \P-a.s.
$$
where $B(r),C_j(r)$ denote the multiplication operators by the functions
$ b'(r,\cdot,X_r(\cdot),u_r)$ and
$\sigma_j'(r,\cdot,X_r(\cdot),u_r)$ respectively. Applying $J_n$ to both sides
it is not hard to conclude that $Y_t^{\epsilon,n}$ has
the Ito differential
$$
dY_s^{\epsilon,n}= A_n\,Y_s^{\epsilon,n}\,ds+ [J_nB(s)Y^\epsilon_s+ J_n\delta^\epsilon b(s)]\,ds
+ [J_nC_j(s)Y^\epsilon_s
    + J_n\delta^\epsilon \sigma_j (s)]\,dW^j_s.
$$
In the following for $y,z\in L^2$, we denote by $y\otimes z$ the rank-one operator
$f\mapsto \<f,z\>_2\,y$ on $L^2$. Using this notation we will consider
 the $\call_2(L^2)$-valued process $Y_s^{\epsilon,n}\otimes Y_s^{\epsilon,n}$, $s\in [t_0, T]$ (recall that if $K$ is a separable Hilbert space, $\call_2(K)$ is the Hilbert space of all bounded linear operators in $X$ for which
 $||X||_{\call_2(K)}^2= tr (X^*X)$ is finite naturally endowed with the product $\<X_1,X_2\>_{\call_2(K)}=tr (X_1^*X_2)$).

 By the Ito formula for Hilbert-space valued Ito processes we have
 $$
 \begin{array}{l}
d(Y_s^{\epsilon,n}\otimes Y_s^{\epsilon,n})= A_n\,(Y_s^{\epsilon,n}\otimes Y_s^{\epsilon,n})\,ds+
 (Y_s^{\epsilon,n}\otimes Y_s^{\epsilon,n})\,A_n^*\,ds
 \\\dis
 \qquad
 +Y_s^{\epsilon,n}\otimes [J_nB(s)Y^\epsilon_s+ J_n\delta^\epsilon b(s)]\,ds
 +[J_nB(s)Y^\epsilon_s+ J_n\delta^\epsilon b(s)]\otimes Y_s^{\epsilon,n}\,ds
 \\\dis
 \qquad
+ Y_s^{\epsilon,n}\otimes[J_nC_j(s)Y^\epsilon_s
    +J_n \delta^\epsilon \sigma_j (s)]\,dW^j_s
    +
  [J_nC_j(s)Y^\epsilon_s
    +J_n \delta^\epsilon \sigma_j (s)]\otimes Y_s^{\epsilon,n}\,dW^j_s
   \\\dis
 \qquad
    +[J_nC_j(s)Y^\epsilon_s
    +J_n \delta^\epsilon \sigma_j (s)]\otimes  [J_nC_j(s)Y^\epsilon_s
    +J_n \delta^\epsilon \sigma_j (s)]\,ds
    \end{array}
$$
and it follows that
 $$
 \begin{array}{l}
Y_s^{\epsilon,n}\otimes Y_s^{\epsilon,n}
\\\dis =
\int_{t_0}^se^{(s-r)A_n}
 \{Y_r^{\epsilon,n}\otimes [J_nB(r)Y^\epsilon_r+ J_n\delta^\epsilon b(r)]
 +[J_nB(r)Y^\epsilon_r+ J_n\delta^\epsilon b(r)]\otimes Y_r^{\epsilon,n}\}
 e^{(s-r)A_n^*}\,dr
 \\\dis
+ \int_{t_0}^se^{(s-r)A_n}\{
Y_r^{\epsilon,n}\otimes[J_nC_j(r)Y^\epsilon_r
    +J_n \delta^\epsilon \sigma_j (r)]
    +
  [J_nC_j(r)Y^\epsilon_r
    +J_n \delta^\epsilon \sigma_j (r)]\otimes Y_r^{\epsilon,n}
    \}
 e^{(s-r)A_n^*}\,dW^j_r
   \\\dis
+\int_{t_0}^se^{(s-r)A_n}
 \{
    [J_nC_j(r)Y^\epsilon_r
    +J_n \delta^\epsilon \sigma_j (r)]\otimes  [J_nC_j(r)Y^\epsilon_r
    +J_n \delta^\epsilon \sigma_j (r)]\}e^{(s-r)A_n^*}\,dr
    \end{array}
$$
The reason for introducing the process $Y^{\epsilon,n}\otimes Y^{\epsilon,n} $
is that we can now make the following computation:
denoting by $tr$ the  trace of operators in $L^2$ we have
$$\E  \< P_{t_0}^N\, Y_{t_0+\epsilon}^{\epsilon,n},
Y_{t_0+\epsilon}^{\epsilon,n}\>_2=
\E\,tr[  P_{t_0}^N\, (Y_{t_0+\epsilon}^{\epsilon,n}\otimes
Y_{t_0+\epsilon}^{\epsilon,n})]
$$
and we can replace $(Y_{t_0+\epsilon}^{\epsilon,n}\otimes
Y_{t_0+\epsilon}^{\epsilon,n})$ by the previous formula. Taking conditional
expectation with respect to $\calf_{t_0}$ the stochastic integral
disappears and we obtain
$$
 \begin{array}{l}
\E  \< P_{t_0}^N\, Y_{t_0+\epsilon}^{\epsilon,n},
Y_{t_0+\epsilon}^{\epsilon,n}\>_2
\\\dis
=\int_{t_0}^{t_0+\epsilon} \E\,tr\Big[  P_{t_0}^N\, e^{(t_0+\epsilon-r)A_n}
 \{Y_r^{\epsilon,n}\otimes [J_nB(r)Y^\epsilon_r+ J_n\delta^\epsilon b(r)]
 \\\dis
 +[J_nB(r)Y^\epsilon_r+ J_n\delta^\epsilon b(r)]\otimes Y_r^{\epsilon,n}\}
 e^{(t_0+\epsilon-r)A_n^*}\Big]\,dr
   \\\dis
+\int_{t_0}^{t_0+\epsilon}\E\,tr\Big[  P_{t_0}^N\,e^{(t_0+\epsilon-r)A_n}
 \{
    [J_nC_j(r)Y^\epsilon_r
    +J_n \delta^\epsilon \sigma_j (r)]\otimes  [J_nC_j(r)Y^\epsilon_r
    +J_n \delta^\epsilon \sigma_j (r)]\}e^{(t_0+\epsilon-r)A_n^*}\Big]\,dr
    \\\dis =2
\E\int_{t_0}^{t_0+\epsilon}   \< P_{t_0}^N\, e^{(t_0+\epsilon-r)A_n}
 [J_nB(r)Y^\epsilon_r+ J_n\delta^\epsilon b(r)],
 e^{(t_0+\epsilon-r)A_n} Y_r^{\epsilon,n}\>_2\,dr
   \\\dis
+ \E\int_{t_0}^{t_0+\epsilon}\< P_{t_0}^N\,e^{(t_0+\epsilon-r)A_n}
    [J_nC_j(r)Y^\epsilon_r
    +J_n \delta^\epsilon \sigma_j (r)], e^{(t_0+\epsilon-r)A_n} [J_nC_j(r)Y^\epsilon_r
    +J_n \delta^\epsilon \sigma_j (r)]\>_2\,dr.
    \end{array}
$$
Next we let $n\to\infty$ and we use the fact that $\|e^{tA_n}f-e^{tA}f\|_2\to 0$ and
$\|J_nf-f\|_2\to 0$ for $f\in L^2$. It follows that
$$
 \begin{array}{l}\dis
\E  \< P_{t_0}^N\, Y_{t_0+\epsilon}^{\epsilon },
Y_{t_0+\epsilon}^{\epsilon}\>_2
=2
\E\int_{t_0}^{t_0+\epsilon}   \< P_{t_0}^N\, e^{(t_0+\epsilon-r)A}
 [ B(r)Y^\epsilon_r+  \delta^\epsilon b(r)],
 e^{(t_0+\epsilon-r)A } Y_r^{\epsilon }\>_2\,dr
   \\\dis
+ \E\int_{t_0}^{t_0+\epsilon}\< P_{t_0}^N\,e^{(t_0+\epsilon-r)A }
    [ C_j(r)Y^\epsilon_r
    +  \delta^\epsilon \sigma_j( r)], e^{(t_0+\epsilon-r)A } [C_j(r)Y^\epsilon_r
    +  \delta^\epsilon \sigma_j (r)]\>_2\,dr.
    \end{array}
$$
Recalling (\ref{cambionorme}), this formula can be written
$$
 \begin{array}{l}\dis
\E  \< P_{t_0}\Pi^N Y_{t_0+\epsilon}^{\epsilon },\Pi^N
Y_{t_0+\epsilon}^{\epsilon}\>
=2
\E\int_{t_0}^{t_0+\epsilon}   \< P_{t_0}\Pi^N e^{(t_0+\epsilon-r)A}
 [ B(r)Y^\epsilon_r+  \delta^\epsilon b(r)],
\Pi^N e^{(t_0+\epsilon-r)A } Y_r^{\epsilon }\>\,dr
   \\\dis
+ \E\int_{t_0}^{t_0+\epsilon}\< P_{t_0}\Pi^N e^{(t_0+\epsilon-r)A }
    [ C_j(r)Y^\epsilon_r
    +  \delta^\epsilon \sigma_j( r)], \Pi^Ne^{(t_0+\epsilon-r)A } [C_j(r)Y^\epsilon_r
    +  \delta^\epsilon \sigma_j (r)]\>\,dr.
    \end{array}
$$
We let $N\to\infty$ and we finally obtain
$$
 \begin{array}{l}\dis
\E  \< P_{t_0} \, Y_{t_0+\epsilon}^{\epsilon },
Y_{t_0+\epsilon}^{\epsilon}\>
=2
\E\int_{t_0}^{t_0+\epsilon}   \< P_{t_0} \, e^{(t_0+\epsilon-r)A}
 [ B(r)Y^\epsilon_r+  \delta^\epsilon b(r)],
 e^{(t_0+\epsilon-r)A } Y_r^{\epsilon }\> \,dr
   \\\dis
+ \E\int_{t_0}^{t_0+\epsilon}\< P_{t_0} \,e^{(t_0+\epsilon-r)A }
    [ C_j(r)Y^\epsilon_r
    +  \delta^\epsilon \sigma_j (r)], e^{(t_0+\epsilon-r)A } [C_j(r)Y^\epsilon_r
    +  \delta^\epsilon \sigma_j (r)]\>\,dr.
    \end{array}
$$
Using the estimate in Proposition \ref{lpfirstvariation} it follows that
$$
\E  \< P_{t_0} \, Y_{t_0+\epsilon}^{\epsilon },
Y_{t_0+\epsilon}^{\epsilon}\>
=\E\int_{t_0}^{t_0+\epsilon}\< P_{t_0} \,e^{(t_0+\epsilon-r)A }
       \delta^\epsilon \sigma_j (r) , e^{(t_0+\epsilon-r)A }  \delta^\epsilon \sigma_j (r)\>\,dr
       +o(\epsilon),
$$
and since $\|e^{tA}f-f\|_4\to 0$ as $t\to 0$ for every $f\in L^4$ we also conclude that
$$
\E  \< P_{t_0} \, Y_{t_0+\epsilon}^{\epsilon },
Y_{t_0+\epsilon}^{\epsilon}\>
=\E\int_{t_0}^{t_0+\epsilon}\< P_{t_0} \,
       \delta^\epsilon \sigma_j (r) ,   \delta^\epsilon \sigma_j (r)\>\,dr
       +o(\epsilon).
$$

Therefore,
in order to finish the proof of (\ref{deltapyydue}), it remains to show that
\begin{equation}\label{quasideltapyydue}
    \E\int_{t_0}^{t_0+\epsilon}\< (P_r-P_{t_0}) \,
       \delta^\epsilon \sigma_j (r) ,   \delta^\epsilon \sigma_j (r)\>\,dr
       =o(\epsilon).
\end{equation}
We fix $\eta\in (0,1/4)$. Since we have $\|\delta^\epsilon \sigma_j(s)\|_4\le C_0$, for some constant  $C_0$,
it follows  that
$$
(-A)^{-\eta}\delta^\epsilon \sigma_j(s)\in
 K:=\{ f\in L^4\;:\; f\in D(-A)^\eta,
\|f\|_{D(-A)^\eta}\le C_0\}.
$$
Since $D(-A)^\eta$ is compactly embedded in $L^4$, the set $K$
is a compact, hence bounded, subset of $L^4$.
We have
$$
| \E  \< (P_r-P_{t_0}) \,
       \delta^\epsilon \sigma_j (r) ,   \delta^\epsilon \sigma_j (r)\>|
\le
\E  \sup_{f\in K  }|\< (P_{r}-P_{t_0})\, (-A)^\eta f,
(-A)^\eta f\>|.
$$
 Since $K$ is compact in $L^4$, for every $\delta>0$
 it can be covered by a finite number $N_\delta$
of open balls with radius $\delta$ and centers denoted $f_i^\delta$, $i=1,\ldots,N_\delta$.
Since $D(-A)^\eta$ is dense in $L^4$, we can assume that $f_i^\delta\in D(-A)^\eta$.
Given $f\in K$, let $i$ be such that $\|f-f_i^\delta\|_4<\delta$; then writing
$$\begin{array}{l}
  \< (P_{r}-P_{t_0}) (-A)^\eta f,(-A)^\eta f\>=
  \< (P_{r}-P_{t_0}) (-A)^\eta f_i^\delta,(-A)^\eta f_i^\delta\>
\\\dis -
  \< (P_{r}-P_{t_0}) (-A)^\eta (f-f_i^\delta),(-A)^\eta (f-f_i^\delta)\>
  +2 \< (P_{r}-P_{t_0}) (-A)^\eta f,(-A)^\eta (f-f_i^\delta)\>
\end{array}
$$
and recalling the notation introduced in (\ref{normextension}) we obtain
$$\begin{array}{l}
|  \< (P_{r}-P_{t_0}) (-A)^\eta f,(-A)^\eta f\>|\le
  |\< (P_{r}-P_{t_0}) (-A)^\eta f_i^\delta,(-A)^\eta f_i^\delta\>|
\\\dis +
  |||P_{r}-P_{t_0}|||\,\delta^2
  +2 |||P_{r}-P_{t_0}|||\,\|f\|_4\,\delta.
\end{array}
$$
Since $K$ is bounded in $L^4$, we conclude that
$$\begin{array}{l}\dis
\sup_{f\in K  }|\< (P_{r}-P_{t_0})\, (-A)^\eta f,
(-A)^\eta f\>\le
\sum_{i=1}^{N_\delta}
 |\< (P_{r}-P_{t_0}) (-A)^\eta f_i^\delta,(-A)^\eta f_i^\delta\>|
\\\dis +
2\sup_{t\in[ t_0,r]}|||P_{t}|||\,\delta^2
  +c \sup_{t\in [t_0,r ]}||| P_{t }|||\,\delta,
\end{array}
$$
for some constant $c$.
Taking expectation,
it follows from (\ref{stimaindaetaperptre})
that
$$
| \E  \< (P_r-P_{t_0}) \,
       \delta^\epsilon \sigma_j (r) ,   \delta^\epsilon \sigma_j (r)\>|\le
\sum_{i=1}^{N_\delta}
\E |\< (P_{r}-P_{t_0}) (-A)^\eta f_i^\delta,(-A)^\eta f_i^\delta\>|
+c(T-r)^{-2\eta}
[\delta^2
  + \delta],
$$
for some constant $c$ independent of $\epsilon$ and $\delta$. By
(\ref{continelleperp}) we conclude that
$$
\limsup_{r \downarrow t_0}
| \E  \< (P_r-P_{t_0}) \,
       \delta^\epsilon \sigma_j (r) ,   \delta^\epsilon \sigma_j (r)\>|
\le
 c(T-t_0)^{-2\eta}
[\delta^2
  + \delta ].
$$
Letting $\delta\to 0$ we conclude that the left-hand side is zero,
and
 (\ref{quasideltapyydue}) follows immediately.

\appendix
\section{Stochastic integrals in $L^p$ spaces}

In this appendix we sketch the construction and some basic properties
of stochastic integrals with respect to a finite dimensional
Wiener process, taking values in an $L^p$-space.
The few facts collected below are enough for the present paper.

Let $(W^1_t,\ldots,W^d_t)_{t\ge0}$ be a standard, $d$-dimensional
Wiener process defined in some complete probability space $(\Omega,\calf,\P)$.
We denote by $(\calf_t)_{t\ge0}$ the corresponding natural filtration,
augmented in the usual way, and we denote by $\calp$ the progressive $\sigma$-algebra
on $\Omega\times [0,T]$, where $T>0$ is a given number.
Let
$L^p:=L^p(D,\cald,m)$ be the usual space,
where $m$ is a positive, $\sigma$-finite measure and $p\in [2,\infty)$.
The integrand processes will be functions $H:\Omega\times [0,T]\times D\to \R^d$, which
are assumed to be $\calp\otimes \cald$-measurable.
When $H$ is of special type, i.e. it has components of the form $$
H^j(\omega,t,x)=\sum_{i=1}^Nh^j_i(\omega,t)f^j_i(x)
$$
for $j=1,\ldots, d$, $h^j_i$ bounded $\calp$-measurable,
$f^j_i$ bounded $\cald$-measurable, then the stochastic integral $I_t(x)$
is defined for fixed $x\in D$ by the formula
$I_t(x)=\int_0^tH_s^j(x)\,dW^j_s= f^j_i(x) \int_0^th^j_i(s)\,dW^j_s
$. Using the Burkholder-Davis-Gundy inequalities
for real-valued stochastic integrals, we have for some constant $c_p$ (depending
only on $p$):
$$
\E |I_t(x)|^p\le c_p\E\left( \int_0^t |H_s(x)|^2ds\right)^{p/2}
$$
where $|H_s(x)|^2=\sum_{j=1}^d|H_s^j(x)|^2$. Since $p\ge 2$ we have, by en elementary inequality,
$$
\E |I_t(x)|^p\le c_p\left( \int_0^t (\E|H_s(x)|^p)^{2/p}ds\right)^{p/2}
= c_p\left( \int_0^t \|H_s(x)\|^2_{L^p(\Omega;\R^d)}ds\right)^{p/2}.
$$
Integrating with respect to $m$ we obtain, again by elementary arguments,
$$
\E \|I_t\|_{L^p(D)}^p\le c_p\int_D \left( \int_0^t (\E|H_s(x)|^p)^{2/p}ds\right)^{p/2}m(dx)
\le
c_p \left( \int_0^t \left(\int_D \E|H_s(x)|^p m(dx)\right)^{2/p}ds\right)^{p/2}
$$
which can be written
\begin{equation}\label{isometryiito}
    \E \|I_t\|_{L^p(D)}^p\le
 c_p\left( \int_0^t (\E \|H_s\|^p_{L^p(D;\R^d)})^{2/p}ds\right)^{p/2}
\end{equation}
 or equivalently
$$
    \|I_t\|_{L^p(\Omega\times D)} \le
 c_p^{1/p}\left( \int_0^t   \|H_s\|^2_{L^p(\Omega\times D;\R^d)}ds\right)^{1/2}.
$$
Finally, by standard arguments, the stochastic integral can be extended to the
class of $\calp\otimes \cald$-measurable integrands $H$ for which the right-hand side of
(\ref{isometryiito}) is finite, and the inequality (\ref{isometryiito}) remains true.

We finally note that from (\ref{isometryiito}) and the H\"older inequality it follows that
\begin{equation}\label{isometryiitoter}
    \E \|I_t\|_{L^p(D)}^p\le
 c_p\int_0^t \E \|H_s\|^p_{L^p(D;\R^d)}ds\; t^{(p-2)/2}.
\end{equation}

Now suppose that there exist  regular conditional probabilities
$\P(\cdot|\calf_t)$
given any $\calf_t$ (this holds for instance if the Wiener process is canonically
realized on the space of $\R^d$-valued continuous functions). Then a slight
modification of the previous passages shows the validity of the following
conditional variant of \eqref{isometryiitoter}: for $0\le r\le t$,
\begin{equation}\label{isometryiitotercond}
    \E^{\calf_r} \|\int_r^tH_s^j\,dW^j_s\|_{L^p(D)}^p\le
 c_p\int_r^t \E^{\calf_r} \|H_s\|^p_{L^p(D;\R^d)}ds\; (t-r)^{(p-2)/2}.
\end{equation}
This is used in the proof of
Proposition \ref{solauxiliaryy}.

\end{document}